\documentclass[11pt]{amsart}
\pdfoutput=1
\usepackage{amsmath}
\usepackage{amssymb}
\usepackage{amsthm} 
\usepackage{bbm}
\usepackage{graphicx}
\usepackage{algorithm}
\usepackage{algorithmic}
\usepackage{clrscode}
\usepackage{fullpage}
\usepackage{stmaryrd}
\usepackage{color}
\usepackage{hyperref}
\usepackage{caption}
\usepackage{subcaption}
\usepackage{bm}
\usepackage{color}
\usepackage{enumerate}
\usepackage{arydshln}
\usepackage{hyperref}

\newtheorem{cor}{Corollary}
\newtheorem{thm}{Theorem}
\newtheorem*{thm*}{Theorem}
\newtheorem{lem}{Lemma}

\newtheorem{Def}{Definition}

\def\x{{\mathbf x}}

\def\b{{\mathbf b}}
\def\y{{\mathbf y}}
\def\n{{\mathbf n}}
\def\m{{\mathbf m}}

\newcommand{\polylog}{\mathop{\rm polylog}}

\DeclareMathOperator*{\argmax}{arg\,max}
\newcommand{\Mod}[1]{\ \text{mod}\ #1}

\begin{document}
\pagestyle{plain}
\newenvironment{frcseries}{\fontfamily{frc} \selectfont}{}
\newcommand{\textfrc}[1]{{\frcseries #1}}
\newcommand{\mathfrc}[1]{\text{\textfrc{#1}}}

\title{Fast Phase Retrieval from Local Correlation Measurements}\thanks{M.A. Iwen:  Department of Mathematics and Department of ECE, Michigan State University ({\tt markiwen@math.msu.edu}).  M.A. Iwen was supported in part by NSF DMS-1416752 and NSA H98230-13-1-0275.\\  \indent A. Viswanathan:  Department of Mathematics, Michigan State University ({\tt aditya@math.msu.edu}). \\  \indent Y. Wang: Department of Mathematics, The Hong Kong University of Science and Technology ({\tt yangwang@ust.hk}). Y. Wang was partially supported by NSF DMS-1043032 and AFOSR FA9550-12-1-0455.}
%
\author{
{Mark Iwen} 
\and Aditya Viswanathan
\and Yang Wang
}

\maketitle
\thispagestyle{empty}

\begin{abstract}
We develop a fast phase retrieval method which can utilize a
large class of local phaseless correlation-based measurements in
order to recover a given signal ${\bf x} \in \mathbb{C}^d$ (up to
an unknown global phase) in near-linear $\mathcal{O} \left( d
\log^4 d \right)$-time.  Accompanying theoretical analysis proves
that the proposed algorithm is guaranteed to deterministically
recover all signals ${\bf x}$ satisfying a natural flatness
(i.e., non-sparsity) condition for a particular choice of
deterministic correlation-based measurements.  A randomized
version of these same measurements is then shown to provide
nonuniform probabilistic recovery guarantees for arbitrary
signals ${\bf x} \in \mathbb{C}^d$.  Numerical experiments
demonstrate the method's speed, accuracy, and robustness in
practice -- all code is made publicly available.

In its simplest form, our proposed phase retrieval method employs a modified lifting scheme akin to the one utilized by the well-known {\em PhaseLift} algorithm.  In particular, it interprets quadratic magnitude measurements of $\x$ as linear measurements of a restricted set of lifted variables, $x_i \overline{x_j}$, for $| j - i | < \delta \ll d$.  This leads to a linear system involving a total of $(2 \delta-1)d$ unknown lifted variables, all of which can then be solved for using only $\mathcal{O}(\delta d)$ measurements.  
Once these lifted variables, $x_i \overline{x_j}$ for $| j - i | < \delta \ll d$, have been recovered, a fast angular synchronization method can then be used to propagate the local phase difference information they provide across the entire vector in order to estimate the (relative) phases of every entry of $\x$.  In addition, the lifted variables corresponding to $x_j \overline{x_j} = |x_j|^2$ automatically provide magnitude estimates for each entry, $x_j$, of $\x$.  The proposed phase retrieval method then approximates $\x$ by carefully combining these entry-wise phase and magnitude estimates.

Finally, we conclude by developing an extension of the proposed method to the 
sparse phase retrieval problem; specifically, we demonstrate 
a sublinear-time compressive phase retrieval algorithm which is
guaranteed to recover a given $s$-sparse vector ${\bf x} \in
\mathbbm{C}^d$ with high probability in just $\mathcal{O}(s
\log^5 s \cdot \log d)$-time using only $\mathcal{O}(s \log^4 s
\cdot \log d)$ magnitude measurements.  In doing so we demonstrate the
existence of compressive phase retrieval algorithms with
near-optimal linear-in-sparsity runtime complexities.

\end{abstract}


\section{Introduction} \label{sec:Intro}

We consider the \textit{phase retrieval problem} of recovering a
vector $\x \in \mathbbm{C}^d$, up to an unknown global phase
factor, from squared magnitude measurements 
\begin{equation} 
  \b := \left \vert M \x \right \vert^2 + \n,
  \label{eq:PR_problem} 
\end{equation}
where $\b \in \mathbbm R^D$ is the vector of phaseless
measurements (with $D\geq d$), $M\in \mathbbm{C}^{D \times d}$ is
a measurement matrix, $\n \in \mathbbm R^D$ denotes measurement
noise, and, $|\cdot |^2: \mathbbm{C}^D \rightarrow \mathbbm{R}^D$
computes the componentwise squared magnitude of each vector
entry. Our objective is to design a computationally efficient and
robust recovery method, $\mathcal A_M: \mathbbm{R}^D \rightarrow
\mathbbm{C}^d$, which can approximately recover ${\bf x}$ using
the magnitude measurements $\b$ that result from any member of a
relatively large class of local correlation-based measurement 
matrices $M \in \mathbbm{C}^{D \times d}$.

\begin{figure}[htbp] \centering
    \includegraphics[scale=0.65]{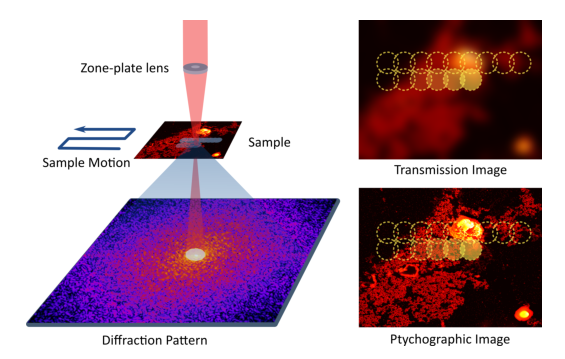}
    \caption{A Typical Imaging Setup for Ptychography.
            \protect\footnotemark}
    \label{fig:ptychography} 
\end{figure}

Phase retrieval problems of this form arise naturally in many crystallography and optics
applications (see, e.g.,
\cite{walther1963question,millane1990phase,harrison1993phase,miao2008extending}).
As an illustrative example, Figure \ref{fig:ptychography}
presents a typical imaging setup for a popular molecular imaging
modality known as ptychography \cite{Rodenburg2008}. 
The figure shows a beam of electromagnetic radiation
illuminating a small region of a test specimen. Under certain
conditions (for example, if the wavelength of the incident
radiation is of the same order as the atomic features in the
specimen), the resulting diffraction pattern contains information
about the atomic structure of the region being imaged.
Therefore, by successively imaging shifts of the specimen, one
might hope to deduce the complete atomic structure from such measurements.  Unfortunately, it can be
shown that the diffraction pattern measured by the detector
corresponds to the (squared) magnitude of the Fourier transform
of the specimen being imaged, meaning, among other things, that all phase information is lost.  This is characteristic of many
molecular imaging methods, where, either due to the underlying
physics or due to instrumentation challenges, the detectors are
not capable of capturing phase information.  Therefore, a phase
retrieval problem needs to be solved in order to recover the
underlying atomic structure.  For ptychography applications, in
particular, and restricting our attention to one dimension, the acquired measurements are of the form
\[  b_\ell(\omega) = \left \vert \mathcal F \left[ \widetilde m ~S_{\ell} x 
        \right] (\omega) \right \vert^2,
\]
where $x$ is the specimen of interest, $\widetilde m$ is a
localized illumination or window function (which depends on the
optical setup), $S_{\ell}$ denotes a shift/translation operator, and $\mathcal F$ denotes the Fourier transform.
Note that the Fourier intensity measurement $b_\ell(\omega)$
corresponds to a specific translate, or shift, $\ell$ of the
unknown specimen $x$.  Furthermore, we may assume that ${\rm
supp}(\widetilde m) \subset {\rm supp}(x)$ since the imaging
field of view is typically much smaller than the support of the
specimen.
\footnotetext{Image credits: Qian, Jianliang, et al. ``Efficient
algorithms for ptychographic phase retrieval.'' Inverse Problems
and Applications, Contemp. Math 615 (2014): 261--280.}

Discretizing the problem on a uniform grid, we obtain 
\begin{equation*} 
  (b_\omega)_\ell = \left \vert \sum_{k=0}^{d-1} \mathbbm 
    e^{\frac{-2 \pi \mathbbm i \, \omega k}d} \left(
    \widetilde m_k \, x_{k+\ell} \right) \right \vert^2, 
    \quad \omega, \ell \in \{0,1,\cdots,d-1\}, 
  \label{eq:ptychography_1d}
\end{equation*}
where $\x, \widetilde \m \in \mathbbm C^d$ and $(b_\omega)_\ell
\in \mathbbm R$ is the analogous discrete measurement
corresponding to a (circular) $\ell$-shift of the unknown vector
$\x$. Due to the finite field of view of the imaging system, we
may also assume that $\widetilde m_k = 0, k>\delta$, where
$\delta \in \mathbbm Z^+$ and $\delta < d$. Hence, we may
write\footnote{All indexing is considered ${\rm mod}\, d$; 
$\overline x$ denotes the complex conjugate of $x$.}
\begin{equation} 
  (b_\omega)_\ell = \left \vert \sum_{k=0}^{\delta-1} \left( 
    \mathbbm e^{\frac{-2 \pi \mathbbm i \, \omega k}d}
    \widetilde m_k \right) x_{k+\ell} \right \vert^2 = 
    \left \vert \sum_{k=0}^{\delta-1} \overline{(\m_\omega)}_k \,
    x_{k+\ell} \right \vert^2, \quad 
    (\m_\omega)_k := \mathbbm e^{\frac{2 \pi \mathbbm 
    i \, \omega k}d} \overline{\widetilde m}_k.
  \label{eq:ptychography_1d_corr_meas}
\end{equation}
The imaging process thus involves collecting measurements $\lbrace
(b_\omega)_\ell\rbrace$ for various Fourier frequencies and 
overlapping image space shifts $(\omega, \ell) \subseteq 
\mathbbm Z^2 \cap [0,d-1]^2$.  Here, in \eqref{eq:ptychography_1d_corr_meas}, and throughout the remainder of the paper, it is important to keep in mind that the parameter $\delta$ always corresponds to the width of each window/mask $\m_\omega$.  The fact that this support size, $\delta$, is significantly smaller that $d$ is what leads to the locality of our proposed measurements.

As we will see in \S \ref{sec:prelims}, our proposed measurement constructions are 
identical to (\ref{eq:ptychography_1d_corr_meas}) for a specific collection of Fourier frequencies and spatial shifts.
We also note that \eqref{eq:ptychography_1d_corr_meas} can be 
interpreted as squared magnitude {\em correlation} 
measurements of the unknown vector $\x$ with the {\em local} (since $\delta \ll d$)
masks $\m_\omega$. In particular, 
\begin{equation}
  \b_\omega = \left \vert {\rm corr}(\m_\omega, \x) 
    \right \vert^2,
    \label{eq:CORRmeas}
\end{equation}
where $\b_\omega \in \mathbbm R^d$ denotes diffraction
measurements corresponding to frequency $\omega$ and all possible 
image space shifts of the unknown vector $\x$.

\subsection{Survey of Previous Work}
\label{subsec:litsurvey}
Given the importance of phase retrieval in crystallographic and
optical imaging methods, there is a rich history of research on
this topic by scientists and practitioners across diverse fields.
With regard to computational algorithms specifically, most
popular methods in use today can trace their origins back to the
alternating projection algorithms developed in the 1970s by
Gerchberg and Saxton \cite{GSaxtonAltProj}, and Fienup
\cite{Fienup_1978_AltProj}. These are iterative formulations
which work by alternately requiring that the following two sets
of constraints be satisfied: 
\begin{enumerate}[{(i)}]
  \item (in data space) the current iterate has the same
    magnitude as that of the measured data, and
  \item (in image space) the current iterate satisfies certain
    problem-specific constraints such as positivity or compact
    support. 
\end{enumerate}
These algorithms are conceptually simple, efficient to implement
(FFT-time for certain measurement constructions) and popular
among the practitioners, despite the lack of a rigorous
mathematical understanding of their properties or global 
recovery guarantees.
Indeed, given the non-convex nature of the phase retrieval
problem, the lack of global convergence results for these methods
is not surprising.  The interested reader is referred to
\cite{Marchesini2007,Fannjiang2012} for some recent developments
and a review of this family of methods, while 
\cite{qian2014efficient,rodenburg2004phase,maiden2009improved}
contain some specific applications of alternating projection
algorithms to Ptychographic imaging.

More recently, there have been significant efforts devoted 
towards 
developing phase retrieval algorithms which are (i)
computationally efficient, (ii) robust to measurement noise, and
(iii) theoretically guaranteed to reconstruct a given vector up
to a global phase error using a near-minimal number of magnitude
measurements.  For example, it has been shown that robust phase
retrieval is possible with $D = \mathcal{O}(d)$ magnitude
measurements by solving a semidefinite programming relaxation
({\em PhaseLift}) of a rank-$1$ matrix recovery problem
\cite{candes2013phaselift,candes2014solving}. This allows
polynomial-time convex optimization methods to be used for phase
retrieval. Furthermore, the runtimes of these convexity-based
methods can be reduced with the use of $\mathcal{O}(d \log d)$
magnitude measurements \cite{demanet2014stable}. Other phase
retrieval approaches include the use of spectral recovery methods
together with magnitude measurement ensembles inspired by
expander graphs \cite{alexeev2014phase}. These methods allow the 
recovery of ${\bf x}$ up to a global phase factor using 
$\mathcal{O}(d)$ noiseless magnitude measurements or $\mathcal
O(d \log d)$ noisy magnitude measurements, and run in 
$\Omega(d^2)$-time in general.\footnote{Their runtime complexity 
is dominated by the time required to solve an overdetermined 
linear system.} Finally, the recently proposed {\em Wirtinger
Flow} algorithm \cite{Candes2015WF} and {\em Truncated Wirtinger 
Flow} (TWF) algorithm \cite{chen2015solving} employ stochastic 
gradient descent schemes with special eigenvector-based 
initialization methods to recover $\x$. The TWF algorithm, for 
example, recovers $\x$ using $\mathcal O(d)$ magnitude 
measurements robustly with computational complexity
\footnote{The
{\em Wirtinger flow} algorithms are known to empirically work
with random masked measurement constructions known as coded
diffraction patterns (CDP) \cite{Candes2014WF,chen2015solving}
in essentially FFT-time, although no robust recovery guarantees 
are available for such measurements.}
$\mathcal O(Dd\log 
1/\epsilon)$, where $\epsilon$ is an accuracy parameter. All 
of these approaches utilize magnitude measurements $\left | M
\x \right |^2$ resulting from either (i) Gaussian random
matrices $M$, (ii) random masked Fourier-based constructions
known as coded diffraction patterns (CDP)
\cite{gross2014improved}, or (iii) unbalanced 
expander graph constructions,
in order to prove their recovery guarantees.

%
\subsection{Main Result}
\label{sec:mainresult}
%
In this paper we demonstrate that a relatively general class of local correlation-based magnitude measurements allow for phase retrieval in
just $\mathcal{O}(d \log^c d)$-time.\footnote{Herein $c$ is a
fixed absolute constant.}  In particular, we construct a
well-conditioned set of Fourier-based measurements which are
theoretically guaranteed to allow for the phase retrieval of a
given vector ${\bf x} \in \mathbb{C}^d$ with high probability in $\mathcal{O}(d \log^4
d)$-time.  These measurements are of particular interest given
that they are closely related to, e.g., ptychography \cite{Rodenburg2008}.  In particular, we prove the following theorem:
\begin{thm*}
({\bf Fast Phase Retrieval from Correlation-Based Measurements})
Let ${\bf x} \in \mathbbm{C}^d$ with $d$ sufficiently large have
$\| \x \|^2_2 \geq C ~ (d~\ln d)^2 \ln^3 (\ln d) ~\| {\bf n}
\|_2$.\footnote{Herein $C, C', C'' \in \mathbbm{R}^+$  are all
fixed and absolute constants.} 
Then, one can select a random
measurement matrix $M \in \mathbbm{C}^{D \times d}$ such
that the following holds with probability at least $1-\frac{1}{C'
\cdot \ln^2(d) \cdot \ln^3 \left( \ln d \right)}$:\footnote{Note that $M$ is selected independently 
of both $\x$ and $\n$.  Furthermore, the nonuniform probability of 
success can be boosted to $1 - p$ for any $p \in (0,1)$ at the expense of introducing additional logarithmic factors of $1/p$ into the 
runtime and measurement bounds.  See \S \ref{sec:FlattenVecs} for additional details.}
The proposed Phase Retrieval Algorithm~\ref{alg:phaseRetrieval} 
(BlockPR) will recover an ${\bf
\tilde{x}} \in \mathbbm{C}^d$ with 
\begin{equation}
\min_{\theta \in [0, 2 \pi)} \left\| {\bf x}  -
\mathbbm{e}^{\mathbbm{i} \theta} {\bf \tilde{x}} \right\|^2_2
~\leq~ C'' (d~\ln d)^2 \ln^3 (\ln d) \| {\bf n} \|_2
\label{MAINthm:ArbRecovery}
\end{equation}
when given arbitrarily noisy input measurements $\b = \left \vert
M \x \right \vert^2 + \n \in \mathbbm{R}^D$ as per
\eqref{eq:PR_problem}.  Here $D$ can be chosen to be
$\mathcal{O}(d \cdot \ln^2(d) \cdot \ln^3 \left( \ln d \right))$.
Furthermore, the proposed Algorithm~\ref{alg:phaseRetrieval} will 
run in $\mathcal{O}(d \cdot \ln^3(d) \cdot \ln^3 \left( \ln d
\right))$-time.
\end{thm*}
To the best of our knowledge, this result provides the best 
existing error guarantee for correlation-based
measurements. Moreover, this result also guarantees exact
recovery, up to a global phase multiple, of $\x$ in the 
noiseless setting (i.e., when $\n = {\bf 0}$). Further, for a particular class of 
flat\footnote{Here, ``flat" simply means that there are no long strings of consecutive entries of $\x$ all of whose magnitudes are less than $\frac{\| \x \|_2}{2 \sqrt{d}}$.  See \S \ref{sec:RecoverFlatVecs} for additional details.} vectors $\x$, $M$ can
be chosen to be a {\em deterministic} matrix arising from local correlation-based measurements, and $D=3d$ 
measurements  suffice for recovery in the noiseless setting. 

Numerical experiments both verify the speed
and accuracy of the proposed phase retrieval approach, as well as
indicate that the approach is highly robust to measurement noise.
Additionally, after establishing and analyzing our general phase
retrieval method, we then utilize it in order to establish the
first known linear-in-sparsity compressive phase retrieval method
capable of recovering $s$-sparse vectors ${\bf x}$ (up to an
unknown phase factor) in only $\mathcal{O}(s \log^c d)$-time.

\subsubsection{An Overview of the Proposed Approach for Local Correlation-Based Measurements}

The proposed phase retrieval approach works in two stages in the ptychographic setting:  During the first lifting (see, e.g., \cite{candes2013phaselift,candes2014solving}) stage, the quadratic magnitude measurements \eqref{eq:CORRmeas} are viewed as linear measurements of $(2 \delta -1)d$ new lifted variables, the entry-wise products $x_i \overline{x_j}$ for $| j - i | < \delta \ll d$.\footnote{Here, $i-j \in \{ -d/2, \dots, -1, 0, 1, \dots, d/2\}$ is always considered modulo $d$.  Note that these particular $(2 \delta -1)d$ entry-wise products are exactly all those that appear when the squared magnitude measurements \eqref{eq:ptychography_1d_corr_meas} are multiplied out for all shifts $\ell$ (treating $\x$ as periodic).}
Given the correlation structure of the original measurements \eqref{eq:CORRmeas}, the resulting lifted linear system in the new lifted unknowns $x_i \overline{x_j}$ is also highly structured.  In particular, the resulting coefficient matrix of the lifted linear system turns out to be block circulant \cite{T2007Eigenvectors}, which both allows explicit condition number bounds to be derived for particular choices of windows $\m_\omega$, and also allows its fast numerical inversion using Fourier transforms.  Thus, we are generally always able to (rapidly) invert the lifted linear system in order to solve for all local entry-wise products $x_i \overline{x_j}$ with $| j - i | < \delta$.  As an immediate consequence, the magnitude of each entry of $\x$ can also be estimated using, e.g., the products $x_j \overline{x_j} = |x_j|^2$.

During the second angular synchronization (see, e.g., \cite{ang_sync_singer}) stage of the algorithm, the local entry-wise products $x_i \overline{x_j}$, $i \neq j$, are used in order to estimate the relative phases of each individual entry of $\x$.  This is done by noting that the phase of each lifted variable $x_i \overline{x_j}$, $\arg \left( x_i \overline{x_j} \right)$, provides an estimate of the phase angle difference between the corresponding entries' phases $\arg(x_i)$ and $\arg(x_j)$.  Hence, the relative phases between all pairs of entries of $\x$ can be approximated by adding these local relative phase differences, provided by $\arg \left( x_i \overline{x_j} \right)$ for $| j - i | < \delta$, together in appropriate (telescoping) sums.  The only potential difficulty with this simple approach occurs when either $x_i$ or $x_j$ happens to be zero (or, more generally, very small in magnitude).  In this case the phase difference given by $\arg \left( x_i \overline{x_j} \right)$ is unreliable/unusable.  In the worst possible case, when, e.g., $\x$ contains $\geq \delta$ consecutive zero entries in a row, the available local phase differences can not span the corresponding gap in reliable phase difference information, and signal recovery (up to a single global phase ambiguity) becomes impossible.  Much of the analysis of this stage focusses on categorizing and circumventing this fundamental difficulty.

Combining both stages above allows one to obtain deterministic algorithms for recovering all $\x$ that do not contain any long contiguous strings of $\geq \delta$ small entries with local correlation measurements.  This result, which applies in the ptychographic setting, is ultimately developed in \S \ref{sec:RecoverFlatVecs} and stated in Theorem~\ref{thm:RecoverFlatVecs}.  Next, in \S \ref{sec:FlattenVecs}, it is shown that the worst case set of signals can still be recovered with high probability if one is allowed to precede one's correlation measurements with an initial randomized global masking operation.  Combining this result with Theorem~\ref{thm:RecoverFlatVecs} yields the main theorem stated above.  Finally, once phase retrieval of sparse signals has been allowed by the employment of randomized masking in \S \ref{sec:FlattenVecs}, it becomes possible to develop sparse phase retrieval algorithms with near-optimal runtimes in \S \ref{sec:PRsparse}.

The remainder of this paper is organized as follows:  In
Section~\ref{sec:prelims} we establish notation and discuss
important preliminary results.  Next, in Section~\ref{sec:Alg},
we present our general phase retrieval algorithm and discuss it's
runtime complexity.  We then analyze our phase retrieval
algorithm and prove recovery guarantees for specific types of
Fourier-based measurement matrices in
Section~\ref{sec:ErrorAnal}.  In Section~\ref{sec:Eval}, we
empirically evaluate the proposed phase retrieval method for
speed and robustness.  Finally, in Section~\ref{sec:PRsparse}, we
use our general phase retrieval algorithm in order to construct a
sublinear-time compressive phase retrieval method which is
guaranteed to recover sparse vectors (up to an unknown phase
factor) in near-optimal time.  Section~\ref{sec:FutureWork} concludes with several suggestions for future work.

\section{Preliminaries:  Notation and Setup}
\label{sec:prelims}

For any matrix $X \in \mathbbm{C}^{D \times d}$ we will denote
the $j^{\rm th}$ column of $X$ by ${\bf X}_j \in
\mathbbm{C}^{D}$.  The conjugate transpose of a matrix $X \in
\mathbbm{R}^{D \times d}$ will be denoted by $X^{*} \in
\mathbbm{C}^{d \times D}$, and the singular values of any matrix
$X \in \mathbbm{C}^{D \times d}$ will always be ordered as
$\sigma_1(X) \geq \sigma_2(X) \geq \dots \geq
\sigma_{\min(D,d)}(X) \geq 0.$  Also, the condition number of the
matrix $X$ will denoted by $\kappa(X) := \sigma_1(X) /
\sigma_{\min(D,d)}(X)$.  We will use the notation $[n] := \{ 1,
\dots, n\} \subset \mathbbm{N}$ for any $n \in \mathbbm{N}$.
Finally, given any ${\bf x} \in \mathbbm{C}^d$, the vector ${\bf
x}^{\rm~ opt}_s \in \mathbbm{C}^d$ will always denote an optimal
$s$-sparse approximation to ${\bf x}$.  That is, it preserves the
$s$ largest entries in magnitudes of ${\bf x}$ while setting the
rest of the entires to $0$.  Note that ${\bf x}^{\rm~ opt}_s \in
\mathbbm{C}^d$ may not be unique as there can be ties for the
$s^{\rm th}$ largest entry in magnitude.

\subsection{An Illustrative Example}
\label{subsec:example}
Before we write down the setup for the general case, we present a
simple but illustrative example which highlights the structure of
our proposed measurements, and provides a general
overview of the reconstruction algorithm. For simplicity, let us
assume that we are given noiseless measurements in this section so that
\[ \b = \left \vert M \x\right \vert^2, \]
with $\x \in \mathbbm C^4$, $\b \in \mathbbm R^{12}$, and $M \in
\mathbbm C^{12 \times 4}$. Further, let us assume that the
measurement matrix $M$ has a block-circulant structure of the
form
\[ 
  M = 
    \left(    
    \begin{array}{cccc}
        \overline{(\m_1)}_1 & \overline{(\m_1)}_2 & 0 & 0 \\
        0 & \overline{(\m_1)}_1 & \overline{(\m_1)}_2 & 0 \\
        0 & 0 & \overline{(\m_1)}_1 & \overline{(\m_1)}_2 \\
        \overline{(\m_1)}_2 & 0 & 0 & \overline{(\m_1)}_1 
        \vspace{0.05in} \\ \hdashline[2pt/2pt]
        & & & \vspace{-0.1in} \\
        \overline{(\m_2)}_1 & \overline{(\m_2)}_2 & 0 & 0 \\
        0 & \overline{(\m_2)}_1 & \overline{(\m_2)}_2 & 0 \\
        0 & 0 & \overline{(\m_2)}_1 & \overline{(\m_2)}_2 \\
        \overline{(\m_2)}_2 & 0 & 0 & \overline{(\m_2)}_1
        \vspace{0.05in} \\ \hdashline[2pt/2pt]
        & & & \vspace{-0.1in} \\
        \overline{(\m_3)}_1 & \overline{(\m_3)}_2 & 0 & 0 \\
        0 & \overline{(\m_3)}_1 & \overline{(\m_3)}_2 & 0 \\
        0 & 0 & \overline{(\m_3)}_1 & \overline{(\m_3)}_2 \\
        \overline{(\m_3)}_2 & 0 & 0 & \overline{(\m_3)}_1 \\
    \end{array}     \right) = 
      \begin{pmatrix}
        M_1 \\ M_2 \\ M_3
    \end{pmatrix},   \] 
where $M_1, M_2, M_3 \in \mathbbm C^{4\times 4}$ are circulant 
matrices and $\m_1, \m_2, \m_3 \in \mathbbm C^4$ are {\em masks} 
or {\em window functions} with finite support. In
particular\footnote{The notation $(m_\ell)_i$ denotes the
$i$-th entry of the $\ell$-th mask.}, 
$(\m_\ell)_i = 0$ for $i=3,4$ and $\ell \in \{1,2,3\}$. The
astute reader will note that this construction describes {\em
correlation} measurements of the unknown vector $\x$ with {\em
local} masks $\m_\ell$; i.e., 
$\b = \left \vert {\rm corr}(\m_\ell, \x) \right \vert^2, \, 
\ell\in \{1,2,3\}$.
Writing out the correlation sum explicitly and setting $\delta =
2$, we obtain 
\begin{equation}
    \left( b_\ell \right)_i = 
        \left \vert \sum_{k=1}^\delta 
        \overline{({\bf m_\ell})}_k \cdot x_{i+k-1} \right \vert^2, 
        \qquad (i, \ell) \in \{1,2,3,4\} \times \{1,2,3\}.
    \label{eq:measurements_toyprob}
\end{equation}
For a suitable choice of mask such as
\begin{equation*}
    (\m_\ell)_k = \left\{ 
    \begin{array}{ll}
        \frac{\mathbbm{e}^{-k/a}}{\sqrt[4]{2\delta -1}} \cdot
            \mathbbm e^{\frac{2 \pi \mathbbm i \cdot (k-1) 
            \cdot (\ell-1)}{2\delta - 1}} & 
                \textrm{if}~k \leq \delta \\ 
        0 &
                \textrm{if}~k > \delta
    \end{array} \right., ~{\rm with}~a:= \max \left\{ 4,~\frac{\delta - 1}{2} \right\},
\end{equation*}
we see that (\ref{eq:measurements_toyprob}) is of the form of 
\eqref{eq:ptychography_1d_corr_meas} in \S \ref{sec:Intro} which described ptychographic measurements.  These are exactly the measurements analyzed in \S \ref{sec:CondNumBounds}.

Having summarized our measurement construction, we now turn our
attention to describing the reconstruction algorithm, which can
be conceptually divided into the following two steps:
\begin{enumerate}[(i)]
  \item Estimate (scaled) local phase differences of the form 
    $\{\overline x_j x_k\, | \, j,k \in [d], |j-k \Mod d|<\delta\}$.
  \item Recover the phases of each entry of $\x$, and consequently $\x$ itself, by using the 
    estimates of the local phase differences from (i).
\end{enumerate}

To obtain the phase differences in step (i) above, we start by rewriting the 
squared correlation measurements (\ref{eq:measurements_toyprob}) 
as follows:
\begin{equation*}
    \left( b_\ell \right)_i = 
        \left \vert \sum_{k=1}^\delta 
            \overline{({\bf m_\ell})}_k \cdot x_{i+k-1} 
            \right \vert^2 
        = \sum_{j,k=1}^\delta ({\bf m_\ell})_j \, \overline{
            ({\bf m_\ell})}_k \, \overline x_{i+j-1} \, 
            x_{i+k-1} 
       := \sum_{j,k=1}^\delta ({\bf m_\ell})_{j,k} \,
            \overline x_{i+j-1} \, x_{i+k-1},
\end{equation*}
where we have used the notation 
$({\bf m_\ell})_{j,k}:=({\bf m_\ell})_j\overline{({\bf m_\ell})}_k$.
This is a linear system of equations for $D=(2\delta-1)d=12$ 
phase differences, $\y \in \mathbbm C^{12}$, with 
\[  {\bf y} =  
\begin{bmatrix}
  |x_1|^2 & \overline x_1x_2 &   
    \overline{x}_2x_1 & |x_2|^2 & \overline{x}_2x_3 & 
    \overline{x}_3x_2 & |x_3|^2 & \overline{x}_3x_4 &
    \overline{x}_4x_3 & |x_4|^2 & \overline{x}_4x_1 & 
        \overline x_1x_4
\end{bmatrix}^T.   \]
By defining $\widetilde \b \in \mathbbm R^{12}$ to be the
interleaved vector of measurements (obtained by permuting the entries
of $\b$)
\[  \widetilde \b = 
\begin{bmatrix}
    (b_1)_1 & (b_2)_1 & (b_3)_1 & 
    (b_1)_2 & (b_2)_2 & (b_3)_2 & 
    (b_1)_3 & (b_2)_3 & (b_3)_3 &
    (b_1)_4 & (b_2)_4 & (b_3)_4
\end{bmatrix}^T,      \]
we may write the resulting linear system as $M^\prime \y =
\widetilde \b$, where
\[ 
  \arraycolsep=1.4pt\def\arraystretch{0.3}
  \! \! M^\prime \! \! \! = \! \! \!
    \left(  \! \! 
    \begin{array}{ccc;{2pt/2pt}ccc;{2pt/2pt}ccc;{2pt/2pt}ccc}
        ({\bf m}_1)_{\mbox{\tiny 1,1}} & 
        ({\bf m_1})_{\mbox{\tiny 1,2}} & 
        ({\bf m}_1)_{\mbox{\tiny 2,1}} & 
        ({\bf m}_1)_{\mbox{\tiny 2,2}} & 
        0 & 0 & 0 & 0 & 0 & 0 & 0 & 0\\
        ({\bf m}_2)_{\mbox{\tiny 1,1}} & 
        ({\bf m}_2)_{\mbox{\tiny 1,2}} & 
        ({\bf m}_2)_{\mbox{\tiny 2,1}} & 
        ({\bf m}_2)_{\mbox{\tiny 2,2}} & 
        0 & 0 & 0 & 0 & 0 & 0 & 0 & 0\\
        ({\bf m}_3)_{\mbox{\tiny 1,1}} & 
        ({\bf m}_3)_{\mbox{\tiny 1,2}} & 
        ({\bf m}_3)_{\mbox{\tiny 2,1}} & 
        ({\bf m}_3)_{\mbox{\tiny 2,2}} & 
        0 & 0 & 0 & 0 & 0 & 0 & 0 & 0
        \vspace{0.05in} \\ \hdashline[2pt/2pt]
        & & & \vspace{-0.02in} \\
        0 & 0 & 0 & 
        ({\bf m}_1)_{\mbox{\tiny 1,1}} & 
        ({\bf m_1})_{\mbox{\tiny 1,2}} & 
        ({\bf m}_1)_{\mbox{\tiny 2,1}} & 
        ({\bf m}_1)_{\mbox{\tiny 2,2}} & 
        0 & 0 & 0 & 0 & 0\\
        0 & 0 & 0 & 
        ({\bf m}_2)_{\mbox{\tiny 1,1}} & 
        ({\bf m}_2)_{\mbox{\tiny 1,2}} & 
        ({\bf m}_2)_{\mbox{\tiny 2,1}} & 
        ({\bf m}_2)_{\mbox{\tiny 2,2}} & 
        0 & 0 & 0 & 0 & 0\\
        0 & 0 & 0 & 
        ({\bf m}_3)_{\mbox{\tiny 1,1}} & 
        ({\bf m}_3)_{\mbox{\tiny 1,2}} & 
        ({\bf m}_3)_{\mbox{\tiny 2,1}} & 
        ({\bf m}_3)_{\mbox{\tiny 2,2}} & 
        0 & 0 & 0 & 0 & 0
        \vspace{0.05in} \\ \hdashline[2pt/2pt]
        & & & \vspace{-0.02in} \\
        0 & 0 & 0 & 0 & 0 & 0 & 
        ({\bf m}_1)_{\mbox{\tiny 1,1}} & 
        ({\bf m_1})_{\mbox{\tiny 1,2}} & 
        ({\bf m}_1)_{\mbox{\tiny 2,1}} & 
        ({\bf m}_1)_{\mbox{\tiny 2,2}} & 
        0 & 0 \\
        0 & 0 & 0 & 0 & 0 & 0 & 
        ({\bf m}_2)_{\mbox{\tiny 1,1}} & 
        ({\bf m}_2)_{\mbox{\tiny 1,2}} & 
        ({\bf m}_2)_{\mbox{\tiny 2,1}} & 
        ({\bf m}_2)_{\mbox{\tiny 2,2}} & 
        0 & 0 \\
        0 & 0 & 0 & 0 & 0 & 0 & 
        ({\bf m}_3)_{\mbox{\tiny 1,1}} & 
        ({\bf m}_3)_{\mbox{\tiny 1,2}} & 
        ({\bf m}_3)_{\mbox{\tiny 2,1}} & 
        ({\bf m}_3)_{\mbox{\tiny 2,2}} & 
        0 & 0
        \vspace{0.05in} \\ \hdashline[2pt/2pt]
        & & & \vspace{-0.02in} \\
        ({\bf m}_1)_{\mbox{\tiny 2,2}} & 
        0 & 0 & 0 & 0 & 0 & 0 & 0 & 0 & 
        ({\bf m}_1)_{\mbox{\tiny 1,1}} & 
        ({\bf m_1})_{\mbox{\tiny 1,2}} & 
        ({\bf m_1})_{\mbox{\tiny 2,1}} \\
        ({\bf m}_2)_{\mbox{\tiny 2,2}} & 
        0 & 0 & 0 & 0 & 0 & 0 & 0 & 0 & 
        ({\bf m}_2)_{\mbox{\tiny 1,1}} & 
        ({\bf m_2})_{\mbox{\tiny 1,2}} & 
        ({\bf m_2})_{\mbox{\tiny 2,1}} \\
        ({\bf m}_3)_{\mbox{\tiny 2,2}} & 
        0 & 0 & 0 & 0 & 0 & 0 & 0 & 0 & 
        ({\bf m}_3)_{\mbox{\tiny 1,1}} & 
        ({\bf m_3})_{\mbox{\tiny 1,2}} & 
        ({\bf m_3})_{\mbox{\tiny 2,1}}
    \end{array}  \! \!   \right) \! \!.    \] 
We draw attention to the block-circulant structure of $M^\prime$,
composed of two blocks $M_1^\prime,M_2^\prime\in \mathbbm C^{3 \times 3}$
highlighted using dashed lines. Much of the speed and elegance of the
proposed algorithm arises from this block-circulant structure:
not only can such systems be inverted in essentially linear-time,
they also lend themselves to easy analysis. As we will see in
Section \ref{sec:CondNumBounds}, we can explicitly write out
condition number bounds for $M^\prime$ resulting from certain
classes of measurement masks.

Note that by solving this linear system, we automatically recover 
the magnitude of $\x$.  In particular, 
\[  |x_1|^2 = y_1, \quad |x_2|^2 = y_4, \quad 
    |x_3|^2 = y_7, \quad |x_4|^2 = y_{10}.  \]
Additionally, after normalizing (to unit magnitude) the entries of 
$\y$, we also recover phase difference estimates,
\[ \phi_{i,j} := \arg(x_j) - \arg(x_i), \quad i,j \in 
    \{1,2,3,4\}, \, |i-j \Mod 4|= 1.  \]
For example, $\phi_{1,2} = \arg(x_2) - \arg(x_1) =
\arg(y_2/|y_2|)$.
We can now recover $\arg(\x)$ as required by step (ii) of our two step recovery algorithm by using a greedy procedure. 

Assume,
without loss of generality, that $|x_1| \geq |x_i|, \, i \in
\{2,3,4\}$. We start by setting $\arg(x_1) = 0$.\footnote{Recall that we can
only recover $\x$ up to an unknown global phase factor which, in this case, will be the true phase of $x_1$.} 
We may now set the phase of $x_2$ and $x_4$ using the estimated
phase differences $\phi_{1,2}$ and $\phi_{1,4}$ respectively;
i.e., 
\[  \arg(x_2) = \arg(x_1) + \phi_{1,2}, \quad 
    \arg(x_4) = \arg(x_1) + \phi_{1,4}.     \] 
Similarly, we next set $\arg(x_3) = \arg(x_2) + \phi_{2,3}$,
thereby recovering all of the entries' unknown phases. We note that the
computational cost of this procedure is essentially linear in the
problem size $d$. Section \ref{sec:Alg} contains a more detailed
description of the algorithm, while Section \ref{sec:ErrorAnal}
includes theoretical recovery guarantees.

\subsection{General Problem Setup}
\label{subsec:setup}
Hereafter we will assume that our measurement matrix $M \in \mathbbm{C}^{D \times d}$ has $D := (2\delta -1) d$ rows corresponding to local correlation-based measurements, where $\delta \in \mathbbm{N}$ represents the support size of the associated correlation masks.  Furthermore, we will utilize the decomposition of $M$ into its $(2\delta -1)$ circulant blocks, $M_1, \dots, M_{2\delta -1} \in \mathbbm{C}^{d \times d}$, given by 
\begin{equation}
M =  \left( \begin{array}{l}  M_1\\ M_2\\ \vdots \\ M_{2\delta - 1}\end{array} \right).
\label{equ:MblockDecomp}
\end{equation}
Here each $M_l \in \mathbbm{C}^{d \times d}$ will be both circulant, with 
\begin{equation}
(M_l)_{i,j} := (\m_l)_{(j-i)~{\rm mod}~d ~+~ 1}
\label{equ:SubMatDef}
\end{equation}
for a mask $\m_l \in \mathbbm{C}^d$, and banded, so that $(\m_l)_{i} = 0$ for all $i > \delta$, and $1 \le l \le 2\delta - 1$.\footnote{Every integer modulo $d$ is considered to be an element of $\{ 0, \dots, d-1\}$ in \S\ref{subsec:setup}.  Furthermore, all indexes of vectors in $\mathbbm{C}^d$ will be considered modulo $d$, $+~1$, as per \eqref{equ:SubMatDef} for the remainder of \S\ref{subsec:setup}.}

As a consequence of this structure, the squared magnitude measurements from the $l^{\rm th}$-block, $\left |  M_l \x \right|^2 \in \mathbbm{R}^d$, can be rewritten as 
\begin{equation}
    \left( |M_l \x|^2 \right)_i = \left(M_l \x \right)_i 
        \overline{ \left( M_l \x \right)_i} = 
        \sum^{\delta}_{j,k = 1} (\m_l)_{j} \overline{(\m_l)}_{k} 
        ~ \overline x_{j+i-1} x_{k+i-1}.
    \label{equ:OrigMagMeas}
\end{equation}
Let $\y \in \mathbbm{C}^{D}$ be defined by
\begin{equation}
    y_{i} :=  \overline x_{\left \lceil \frac{i+\delta-1}
        {2\delta-1} \right \rceil} 
        x_{\left \lceil \frac{i+\delta-1}{2\delta-1} \right 
        \rceil + ((i+\delta-2)~{\rm mod}~(2\delta-1)) - \delta+1}.
    \label{equ:Defy}
\end{equation}
Furthermore, let ${\bf 0}_\alpha \in \mathbbm{R}^{1 \times \alpha}$ be the row vector of $\alpha$ zeros for any given $\alpha \in \mathbbm{N}$, and let $\tilde{{\bf m}}_{(l,j)} \in \mathbbm{C}^{1 \times \delta}$ be such that
\begin{equation}
\left( \tilde{{\bf m}}_{(l,j)} \right)_k :=  (\m_l)_{j} \overline{(\m_l)}_{k}.
\end{equation}
We can now re-express $|  M_l \x |^2 \in \mathbbm{R}^d$ from \eqref{equ:OrigMagMeas} as $\tilde{M}_l \y$, where $\tilde{M}_l \in \mathbbm{C}^{d \times D}$ is a $(2\delta-1)$-circulant matrix defined by
\begin{equation*}
\left( \begin{array}{lllllllllll}  \tilde{{\bf m}}_{(l,1)} & {\bf 0}_{\delta-2} & \tilde{{\bf m}}_{(l,2)}  & {\bf 0}_{\delta-2} & \tilde{{\bf m}}_{(l,3)} & \hdots & \tilde{{\bf m}}_{(l,\delta)} & 0 & 0 & \hdots & 0\\
{\bf 0}_{2\delta-1} & \tilde{{\bf m}}_{(l,1)} & {\bf 0}_{\delta-2} & \tilde{{\bf m}}_{(l,2)}  & {\bf 0}_{\delta-2} & \tilde{{\bf m}}_{(l,3)} & \hdots & \tilde{{\bf m}}_{(l,\delta)} & 0 & \hdots & 0\\
&  &  &  & & \ddots & & &  &  & \\
\left( \tilde{{\bf m}}_{(l,2)} \right)_2 & \hdots & \left( \tilde{{\bf m}}_{(l,2)} \right)_\delta &  {\bf 0}_{\delta-2} & \tilde{{\bf m}}_{(l,3)} & {\bf 0}_{\delta-2}  & \hdots & 0 &  \tilde{{\bf m}}_{(l,1)} & {\bf 0}_{\delta-2} & \left( \tilde{{\bf m}}_{(l,2)} \right)_1 \\
\end{array} \right).
\end{equation*}

Finally, after reordering the entries of $\left |  M \x \right|^2$ via a permutation matrix $P \in \{ 0,1 \}^{D \times D}$, we arrive at our final form
\begin{equation}
P \left |  M \x \right|^2 = M' \y = \left( \begin{array}{llllllll}  M'_1& M'_2 & \hdots & M'_{\delta} & 0 & 0 & \hdots & 0\\
0 & M'_1& M'_2 & \hdots & M'_{\delta} & 0 & \hdots & 0\\ 
 & &  & \ddots &  &  &  & \\
 M'_2  & \hdots & M'_{\delta} & 0 & \hdots & 0 & \hdots & M'_1  
\end{array} \right) \y. 
\label{equ:LinProb}
\end{equation}
Here $M' \in \mathbbm{C}^{D \times D}$ is a block circulant matrix \cite{T2007Eigenvectors} whose blocks, $M'_1, \dots, M'_{\delta} \in \mathbbm{C}^{(2\delta-1) \times (2\delta-1)}$, have entries
\begin{equation}
(M'_l)_{i,j} := \left\{ \begin{array}{ll} (\m_i)_{l} \overline{(\m_i)}_{j+l-1} & \textrm{if}~ 1 \leq j \leq \delta - l +1 \\ 0 & \textrm{if}~ \delta - l + 2 \leq j \leq 2\delta - l -1 \\  (\m_i)_{l+1} \overline{(\m_i)}_{l + j - 2 \delta + 1} & \textrm{if}~2 \delta - l \leq j \leq 2\delta - 1, ~\textrm{and}~l < \delta \\ 0 & \textrm{if}~ j > 1, ~\textrm{and}~ l = \delta  \end{array} \right..
\label{eq:M'DEF}
\end{equation}

Let $I_\alpha$ denote the $\alpha \times \alpha$ identity matrix.  We now note that $M'$ can be block diagonalized via the unitary block Fourier matrices $U_{\alpha} \in \mathbbm{C}^{\alpha d \times \alpha d}$, with parameter $\alpha \in \mathbbm{N}$, defined by
\begin{equation}
U_{\alpha} := \frac{1}{\sqrt{d}} \left( \begin{array}{llll}  I_\alpha & I_\alpha & \hdots & I_\alpha \\
 I_\alpha & I_\alpha \mathbbm{e}^{\frac{2 \pi \mathbbm{i}}{d}} & \hdots  & I_\alpha \mathbbm{e}^{\frac{2 \pi \mathbbm{i} \cdot (d-1)}{d}}\\ 
 &  & \ddots   & \\
  I_\alpha &  I_\alpha \mathbbm{e}^{\frac{2 \pi \mathbbm{i} \cdot (d-2)}{d}} & \hdots & I_\alpha \mathbbm{e}^{\frac{2 \pi \mathbbm{i} \cdot (d-2) \cdot (d - 1)}{d}}  \\
 I_\alpha &  I_\alpha \mathbbm{e}^{\frac{2 \pi \mathbbm{i} \cdot (d-1)}{d}} & \hdots & I_\alpha \mathbbm{e}^{\frac{2 \pi \mathbbm{i} \cdot (d-1) \cdot (d - 1)}{d}}  
\end{array} \right). 
\label{equ:UDef}
\end{equation}
More precisely, one can see that we have
\begin{equation}
U^*_{2\delta - 1}~ M' ~U_{2\delta - 1} = J :=  \left( \begin{array}{lllll}  J_1& 0 & 0 & \hdots & 0\\
0 & J_2 & 0 & \hdots & 0\\
 & &  \ddots &  &   \\
 0 & 0 & 0 & J_{d-1} & 0\\
0 & 0 & 0 & 0 & J_{d}\\
\end{array} \right)
\label{equ:Diagonalized}
\end{equation}
where $J \in \mathbbm{C}^{D \times D}$ is block diagonal with blocks $J_1, \cdots, J_d \in \mathbbm{C}^{(2\delta - 1) \times (2\delta - 1)}$ given by
\begin{equation}
J_{k} := \sum^{\delta}_{l = 1} M'_{l} \cdot \mathbbm{e}^{\frac{2 \pi
  \mathbbm{i} \cdot (k-1) \cdot (l-1)}{d}}.
\label{def:JkMat}
\end{equation}
Not so surprisingly, the fact that any block circulant matrix can be block diagonalized by block Fourier matrices will lead to more efficient computational techniques below.

\subsection{Johnson-Lindenstrauss Embeddings and Restricted Isometries}

Below we will utilize results concerning \textit{Johnson-Lindenstrauss embeddings} \cite{JLoriginal,frankl1988johnson,achlioptas2001,dasgupta2003elementary,baraniuk2008simple,krahmer2010new} of a given finite set $\mathcal{S} \subset \mathbbm{C}^d$ into $\mathbbm{C}^m$ for $m < d$.  These are defined as follows:
\begin{Def}
Let $\epsilon \in (0, 1)$, and $\mathcal{S} \subset \mathbbm{C}^d$ be finite.  An $m \times d$ matrix $A$ is a linear Johnson-Lindenstrauss embedding of $\mathcal{S}$ into $\mathbbm{C}^m$ if
$$(1 - \epsilon) \| ~{\bf u} - {\bf v}~ \|_2^2 \leq \| ~A{\bf u} - A{\bf v}~ \|_2^2 \leq (1 + \epsilon) \| ~{\bf u} - {\bf v}~ \|_2^2$$
holds $\forall {\bf u},{\bf v} \in \mathcal{S} \cup \{\bf 0\}$.  In this case we will say that $A$ is a JL($m$,$d$,$\epsilon$)-embedding of $\mathcal{S}$ into $\mathbbm{C}^m$.
\label{Def:JL}
\end{Def}  
\noindent Linear JL($m$,$d$,$\epsilon$)-embeddings are closely related to the \textit{Restricted Isometry Property} \cite{candes2005decoding,baraniuk2008simple,HolgerBook}.
\begin{Def}
Let $s \in [d]$ and $\epsilon \in (0, 1)$. The matrix $A \in \mathbbm{C}^{m \times d}$ has the Restricted Isometry Property if 
\begin{equation}
(1 - \epsilon) \| ~{\bf x}~\|_2^2 \leq \| ~A {\bf x}~\|_2^2 \leq (1 + \epsilon) \| ~{\bf x}~ \|_2^2
\label{equ:RIP}
\end{equation}
holds $\forall {\bf x} \in \mathbbm{C}^d$ containing at most $s$ nonzero coordinates.  In this case we will say that $A$ is RIP($s$,$\epsilon$). 
\end{Def}

\noindent In particular, the following theorem due to Krahmer and Ward \cite{krahmer2010new,HolgerBook} demonstrates that a matrix with the restricted isometry property can be used to construct a Johnson-Lindenstrauss embedding matrix.

\begin{thm}
Let $\mathcal{S} \subset \mathbbm{C}^d$ be a finite point set with $|\mathcal{S}| = M$.
For $\epsilon,p \in (0,1)$, let $A \in \mathbbm{C}^{m \times d}$ be RIP($2s$,$\epsilon / C_1$) for some $s \geq C_2 \cdot \ln(4M/p)$.\footnote{Here $C_1,C_2 \in (1,\infty)$ are both fixed absolute constants.}  Finally, let $B \in \left\{ -1,0,1\right\}^{d \times d}$ be a random diagonal matrix with independent and identically distributed (i.i.d.) symmetric Bernoulli entries on its diagonal.  Then, $AB$ is a JL($m$,$d$,$\epsilon$)-embedding of $\mathcal{S}$ into $\mathbbm{C}^m$ with probability at least $1-p$.
\label{thm:KW}
\end{thm}

Below we will utilize Theorem~\ref{thm:KW} together with a result concerning the restricted isometry property for sub-matrices of a Fourier matrix.  Let $F \in \mathbbm{C}^{d \times d}$ be the unitary $d \times d$ discrete Fourier transform matrix.  The \textit{random sampling matrix}, $R' \in \mathbbm{C}^{m \times d}$,  for $F$ is then
\begin{equation}
R'  := \sqrt{\frac{d}{m}} \cdot RF
\label{equ:RandSamp}
\end{equation}
where $R \in \{ 0,1\}^{m \times d}$ is a random matrix with exactly one nonzero entry per row (i.e., each entry's column position is drawn independently from $[d]$ uniformly at random with replacement).  The following theorem is proven in \cite{HolgerBook}.\footnote{See Theorem 12.32 in Chapter 12.}

\begin{thm}
Let $p \in (0,1)$.  If the number of rows in the random sampling matrix $R' \in \mathbbm{C}^{m \times d}$ satisfies both 
\begin{equation}
\frac{m}{\ln(9 m)} \geq C_3 \cdot \frac{s \ln^2(8s) \ln(8d)}{\epsilon^2}
\end{equation}
and
\begin{equation}
m \geq C_4 \cdot \frac{s \log(1/p)}{\epsilon^2},
\end{equation}
then $R'$ will be RIP($2s$,$\epsilon / C_1$) with probability at least $1-p$.\footnote{Here $C_3,C_4 \in (1,\infty)$ are both fixed absolute constants.}
\label{thm:FourierBON}
\end{thm}
We are now prepared to present and analyze our phase retrieval method.


\section{{\em BlockPR}: A Fast Phase Retrieval Algorithm} \label{sec:Alg}

The proposed phase retrieval algorithm ({\em BlockPR}) works in two stages.  In
the first stage, the vector ${\bf y} \in \mathbbm{C}^D$ from
\eqref{equ:Defy} of local entrywise products of ${\bf x} \in
\mathbbm{C}^d$ with its conjugate is approximated by solving the
linear system \eqref{equ:LinProb}.  That is, we compute $(M')^{-1} P \b$
where $\b \in \mathbb{R}^D$ are our noisy local correlation-based measurements from \eqref{eq:PR_problem}, and $M'$ and $P$ are as in \eqref{equ:LinProb}.
This yields
\begin{equation}
{\bf \tilde{y}} := (M')^{-1} P \b ~=~ (M')^{-1} P \left \vert M \x \right \vert^2 + (M')^{-1} P \n ~=~ \y + (M')^{-1} P \n.
\label{eq:InvertMeas}
\end{equation}
where $\y$ is as in \eqref{equ:Defy}.  

Next, a
greedy algorithm is used to recover the magnitudes and phases of
each entry of ${\bf x}$ (up to a global phase
factor) from our estimate of ${\bf y}$.  To see how this works, note that ${\bf y}$ will contain
all of the products $x_{i} \overline{x}_{j}$ for all $i,j \in
[d]$ with $|(i-j) ~{\rm mod}~d |< \delta$.\footnote{In \S \ref{sec:Alg} and below we always consider any integer modulo $d$ to be in the set $\{ - \lceil \frac{d}{2} \rceil + 1, \dots, \lfloor \frac{d}{2} \rfloor \}$, for $\delta < \lfloor \frac{d}{2} \rfloor$.}  As a result, the
magnitude of each entry $x_j$ can be estimated using the entry of $\y$ corresponding to
$x_{j} \overline{x}_{j} = | x_j |^2$.  Similarly, as long as both
$x_{j} \overline{x}_{j} > 0$ and $x_{i} \overline{x}_{i} > 0$ hold, one can also compute the phase
difference $\arg(x_i) - \arg(x_j)$ from $\arg \left( x_i
  \overline{x}_{j} \right)$.  Thus, the phase of each $x_i$ can be
  determined once $\arg(x_j)$ is established for a neighboring entry.  Repeating this
  process allows one to determine a network of phase differences
  which all depend uniquely on the choice of a single entry's
  unknown phase.  This entry's phase becomes the global phase
  factor $\mathbbm{e}^{\mathbbm{i} \theta}$ from \eqref{MAINthm:ArbRecovery}.  See Algorithm~\ref{alg:phaseRetrieval} for
  additional details.

\begin{algorithm}
  \renewcommand{\algorithmicrequire}{\textbf{Input:}}
  \renewcommand{\algorithmicensure}{\textbf{Output:}}
  \caption{ {\em BlockPR}} \label{alg:phaseRetrieval}
  \begin{algorithmic}[1] \REQUIRE Local Correlation Measurements $\b \in
    \mathbbm{R}^D$ (Recall \eqref{eq:PR_problem}) \ENSURE ${\bf \tilde{x}} \in
    \mathbbm{C}^d$ with ${\bf \tilde{x}} \approx
    \mathbbm{e}^{-\mathbbm{i} \theta} {\bf x}$ for some $\theta
    \in [0, 2 \pi]$ 
    \STATE Compute
    ${\bf \tilde{y}} := (M')^{-1} P \b $ (see
    \eqref{eq:InvertMeas}) \STATE Use Algorithm~\ref{alg:angSync}
    with input ${\bf \tilde{y}} \in \mathbbm{C}^D$ to compute $\tilde{\phi}_j \approx \phi_j := \arg(x_j)$ for all $j \in [d]$
    \STATE Set $\tilde{x}_j = \sqrt{|\tilde{y}_{j'}|}
    \mathbbm{e}^{\mathbbm{i} \tilde{\phi}_j}~\forall j \in [d]$, where $\tilde{y}_{j'}$ (computed in line 1) is s.t.
    $\tilde{y}_{j'} = x_{j} \overline{x}_{j} + \left( (M')^{-1} P \n \right)_{j'}$
  \end{algorithmic} \end{algorithm}

Recall from \eqref{eq:InvertMeas} that we have access to
\begin{equation}
{\bf \tilde{y}} ~=~ \y + {\bf \tilde{n}} \in \mathbbm{C}^D.
\label{eq:InvertMeasDefs}
\end{equation}
where ${\bf \tilde{n}} := (M')^{-1} P \n$ results from the measurement noise, and $D = (2 \delta - 1)d$.  Thus, from \eqref{equ:Defy} we have an index $k(i,j) \in [D]$ for all $i,j \in [d]$ with $|(i-j) ~{\rm mod}~d |< \delta$ such that
\begin{equation}
\tilde{y}_{k(i,j)} ~=~ x_{i}\overline{x}_{j} + \tilde{n}_{k(i,j)}.
\label{eq:Defk(i,j)}
\end{equation}
We will utilize the index function $k$ for ${\bf \tilde{y}}$ defined by \eqref{eq:Defk(i,j)} in Algorithm~\ref{alg:angSync}.

\begin{algorithm}
  \renewcommand{\algorithmicrequire}{\textbf{Input:}}
  \renewcommand{\algorithmicensure}{\textbf{Output:}}
  \caption{Greedy Angular Synchronization}
  \label{alg:angSync} 
  \begin{algorithmic}[1] 
    \REQUIRE $\tilde{y}_{k(i,j)} ~=~ x_{i}\overline{x}_{j} + \tilde{n}_{k(i,j)} \quad \forall i,j \in [d] \textrm{~with~} |(i-j) ~{\rm mod}~d |< \delta$.  
    \ENSURE Phase angles: $\tilde{\phi}_j \approx \phi_j := \arg(x_j)$ for all $j \in [d]$ (up to a global phase).  
   
    \STATE \textit{\% Estimate largest magnitude entry and set its phase to zero.}  
    \STATE $a := \argmax_j ~ |\tilde{y}_{k(j,j)}| =  \argmax_j \left(\left| ~|x_j|^2 + \tilde{n}_{k(j,j)} ~\right| \right), \quad j = 0,\dots,d-1.$ 
    \STATE  $\tilde{\phi}_a \leftarrow 0$ ~~{\em \% Note: We approximate the unknown phases up to a global phase factor.}  \vspace{0.1 in}
  
  \STATE \textit{ \%Define a binary vector,
  $\mbox{phaseFlag} \in \lbrace0,1 \rbrace^d$, to keep track of
  entries whose phase has already\\ ~\%been set.  This ensures that each phase is estimated once, and then not changed again.} \[
    \mbox{phaseFlag}_i = \left \lbrace \begin{array}{cc} 0,
    \qquad & i=a, \\ 1, \qquad & \mbox{else}.  \end{array}
\right. \]  
\STATE Initialize $j \leftarrow a$;  \vspace{.1in}
  
\STATE \textit{\% Estimate phases of all entries of $\x$}
\WHILE{ $\displaystyle \sum_{i \in (j,j+ \delta)} \mbox{phaseFlag}_i > 0$ } \vspace{.1 in}

\STATE \textit{\% Estimate phases of $2\delta -1$ entries nearest current $x_j$}
\FOR{ $i = 1-\delta, 2-\delta, \dots,0,\dots,\delta-1$ } \vspace{.1in} 

\STATE \textit{\% Do not over-write previously estimated phases}
\IF{ $\mbox{phaseFlag}_{j+i ~{\rm mod}~d} == 1$ }  \vspace{.1in}

\STATE \textit{\%Use the current
  reference phase, $\tilde{\phi}_j$, and the input phase difference estimates,}\\ \textit{\% $\arg\left( \tilde{y}_{k(j+i ~{\rm mod}~d,~j)} \right) \approx \arg \left( x_{j+i ~{\rm mod}~d}
  \overline{x}_{j} \right)$ and $\arg\left( \tilde{y}_{k(j,~j+i ~{\rm mod}~d)} \right) \approx$}\\
  \textit{\% $\arg \left( x_{j}\overline{x}_{j+i ~{\rm mod}~d} \right)$, to estimate the phase of the entry $x_{j+i ~{\rm
  mod}~d}$.}
  
\STATE  $\tilde{\phi}_{j+i ~{\rm mod}~d} \leftarrow \tilde{\phi}_j + \frac{1}{2} \left( \arg \left( \tilde{y}_{k(j+i ~{\rm mod}~d,~j)} \right) - \arg \left( \tilde{y}_{k(j,~j+i ~{\rm mod}~d)} \right) \right)$;  \vspace{0.1 in}
\STATE \textit{\% Remember that the phase of entry $x_{j+i ~{\rm mod}~d}$ has now been estimated}
\STATE $\mbox{phaseFlag}_{j+i ~{\rm mod}~d} \leftarrow 0$;    \vspace{0.1 in}

\ENDIF 

\ENDFOR 

\STATE \textit{\% Update the reference entry to be the largest neighboring entry of the current $x_j$} \[ j \leftarrow \left( j + \argmax_{0 < i < \delta}~\left| \tilde{y}_{k(j+i ~{\rm
mod}~d,~j+i ~{\rm mod}~d)} \right| \right) ~{\rm mod}~d ;\]

\ENDWHILE \end{algorithmic} \end{algorithm}

It is important to note that Algorithm~\ref{alg:phaseRetrieval}
assumes that the block circulant matrix $M'$ arising from our
choice of measurements, $M$, is invertible.  As we shall see in
\S\ref{sec:ErrorAnal} and \S\ref{sec:Eval}, this is relatively
easy to achieve.  Similarly, Algorithm~\ref{alg:angSync}
implicitly assumes that ${\bf \tilde{y}}$ does not contain any strings of
$\delta - 1$ consecutive zeros (or, more generally, $\delta - 1$
consecutive entires with ``very small'' magnitudes).  This
assumption will also be discussed in \S\ref{sec:ErrorAnal} and
\S\ref{sec:Eval}, and justified for measurements arising from arbitrary ${\bf x}$ by
modifying the measurements $M$.   For the time being, then, we
are left free to consider the computational complexity of
Algorithm~\ref{alg:phaseRetrieval}.

\subsection{Runtime Analysis}
\label{sec:AlgRuntimeAnal}

We will begin our analysis of the runtime complexity of
Algorithm~\ref{alg:phaseRetrieval} by considering the computation
of ${\bf \tilde y} \in \mathbbm{C}^{D}$ in line 1.  Recalling
\S\ref{sec:prelims}, we note that the permutation matrix $P$ is
based on a simple row reordering that clusters the first rows of
$M_1, \dots, M_{2\delta - 1}$ into a contiguous block, the second
rows of $M_1, \dots, M_{2\delta - 1}$ into a second contiguous
block, etc. (see \eqref{equ:MblockDecomp} and
\eqref{equ:SubMatDef}).  Thus, $P\b$ is simple to
compute using only $\mathcal{O}(d \cdot \delta)$-operations.  To
finish calculating ${\bf \tilde{y}} = (M')^{-1} P \b$ we then use the decomposition of $M'$ from
\eqref{equ:Diagonalized} and compute ${\bf \tilde{y}} = U_{2 \delta - 1}
J^{-1} U^*_{2 \delta - 1}  P \b$.  

Recalling the definition of $U_{2 \delta - 1}$ \eqref{equ:UDef},
one can see that both $U_{2 \delta - 1}$ and $U^*_{2 \delta - 1}$
have fast matrix-vector multiplies (i.e., because they can be
computed by performing $2 \delta-1$ independent fast Fourier
transforms on different sub-vectors of size $d$).  Hence,
matrix-vector multiplies with both of these matrices can be
accomplished with $\mathcal{O}( \delta \cdot d \log d)$
operations.  Finally, $J$ is block-diagonal with $d$ blocks of
size $(2 \delta - 1) \times (2 \delta - 1)$ (see
\eqref{def:JkMat}).  Thus, $J$ and $J^{-1}$ can both be computed
using $\mathcal{O}(d \cdot\delta^3)$ total operations.  Putting
everything together, we can now see that line 1 of
Algorithm~\ref{alg:phaseRetrieval} requires only $\mathcal{O}(d
\cdot \delta^3 +  \delta \cdot d \log d)$ operations in general.
Furthermore, these computations can easily benefit from
parallelism due to the fact that the calculations above are all
based on explicitly defined block decompositions.

The second line of Algorithm~\ref{alg:phaseRetrieval} calls
Algorithm~\ref{alg:angSync} whose runtime complexity is dominated
by its main while-loop (lines 7 through 19).  This loop will
visit each entry of the input vector ${\bf y}$ at most a constant
number of times.  Hence, it requires $\mathcal{O}(\delta \cdot
d)$ operations.  Finally, the third line of
Algorithm~\ref{alg:phaseRetrieval} uses only $\mathcal{O}(d)$
operations.  Thus, the total runtime complexity of
Algorithm~\ref{alg:phaseRetrieval} is $\mathcal{O}(d \cdot
\delta^3 +  \delta \cdot d \log d)$ in general.

\section{Error Analysis and Recovery Guarantees}
\label{sec:ErrorAnal}

In this section we analyze the performance of the proposed phase retrieval method (see Algorithm~\ref{alg:phaseRetrieval}), and demonstrate measurement matrices which allow it to recover arbitrary vectors, up to an unknown phase factor, with high probability.  Our analysis proceeds in two steps.  First, in \S\ref{sec:CondNumBounds} and \S\ref{sec:RecoverFlatVecs}, we demonstrate the existence of a deterministic set of correlation-based measurements, $M \in \mathbbm{C}^{D \times d}$, which allow Algorithm~\ref{alg:phaseRetrieval} to recover all relatively flat (i.e., non-sparse) vectors ${\bf x} \in \mathbbm{C}^d$.  Herein, ``flat'' simply means that there are no long strings of consecutive entries of ${\bf x}$ all of whose magnitudes are less than $\frac{\| {\mathbf x}  \|_2}{2\sqrt{d}}$.  The proposed measurements $M$ are Fourier-like, roughly corresponding to a set of damped and windowed Fourier measurements of overlapping portions of ${\bf x}$.  In  addition to being well conditioned, these Fourier measurements also have fast inverse matrix-vector multiplies via (an additional usage of) the FFT.  Hence, they confer additional computational advantages beyond those already enjoyed by our general block circulant measurement setup.

Next, in \S\ref{sec:FlattenVecs}, we extend our deterministic recovery guarantee for flat vectors to a nonuniform probabilistic recovery guarantee for arbitrary vectors.  This is accomplished by right-multiplying $M$ with a concatenation of several Johnson-Lindenstrauss embedding matrices, each of which tends to ``flatten out'' vectors they are multiplied against.  In particular, we construct a set of such matrices which are both (i) collectively unitary, and (ii) rapidly invertible as a group via (yet another usage of) the FFT.  The fact that this flattening matrix is unitary preserves the well conditioned nature of our initial measurements, $M$.  Furthermore, the fact that the flattening matrix enjoys a fast inverse matrix-vector multiply via the FFT allows us to maintain computational efficiency.  Finally, the fact that the flattening matrix produces a flattened version of ${\bf x}$ with high probability allows us to apply our deterministic recovery guarantee for flat vectors to vectors which are not initially flat.  The end result of this line of reasoning is the following recovery guarantee for noisy measurements.

\begin{thm}
Let ${\bf x} \in \mathbbm{C}^d$ with $d$ sufficiently large have $\| \x \|^2_2 \geq C ~ (d~\ln d)^2 \ln^3 (\ln d) ~\| {\bf n} \|_2$.\footnote{Herein $C, C', C'' \in \mathbbm{R}^+$  are all fixed and absolute constants.}  Then, one can select a random measurement matrix $\tilde{M} \in \mathbbm{C}^{D \times d}$ independently of both $\x$ and $\n$ such that the following holds with probability at least $1-\frac{1}{C' \cdot \ln^2(d) \cdot \ln^3 \left( \ln d \right)}$:  Algorithm~\ref{alg:phaseRetrieval} will recover an ${\bf \tilde{x}} \in \mathbbm{C}^d$ with 
\begin{equation}
\min_{\theta \in [0, 2 \pi]} \left\| {\bf x}  - \mathbbm{e}^{\mathbbm{i} \theta} {\bf \tilde{x}} \right\|^2_2 ~\leq~ C'' (d~\ln d)^2 \ln^3 (\ln d) \| {\bf n} \|_2
\label{thm:ArbRecovery}
\end{equation}
when given arbitrarily noisy input measurements $\b = \left \vert M \x \right \vert^2 + \n \in \mathbbm{R}^D$ as per \eqref{eq:PR_problem}.  Here $D$ can be chosen to be $\mathcal{O}(d \cdot \ln^2(d) \cdot \ln^3 \left( \ln d \right))$.  Furthermore, Algorithm~\ref{alg:phaseRetrieval} will run in $\mathcal{O}(d \cdot \ln^3(d) \cdot \ln^3 \left( \ln d \right))$-time.
\label{thm:RecoverArbVecs}
\end{thm}

Note that the error bound in \eqref{thm:ArbRecovery} is probably
sub-optimal due to, e.g., the appearance of the $d^2\polylog d$-term on its
righthand side.  It is certainly the case, at least, that
stronger error bounds are achievable using complex normal random
measurements in combination with optimization techniques (see,
e.g., \cite{candes2014solving,CandesFix}).  None the less,
Theorem~\ref{thm:RecoverArbVecs} provides the best existing error
guarantee the authors are aware of for local correlation-based
measurements, and computational experiments indicate that
Algorithm~\ref{alg:phaseRetrieval} is highly robust to
measurement noise in practice (see \S\ref{sec:Eval}).
Furthermore, it is important to point out that
Theorem~\ref{thm:RecoverArbVecs} also guarantees exact
recovery, up to a global phase multiple, of $\x$ by
Algorithm~\ref{alg:phaseRetrieval} in the noiseless setting
(i.e., when $\n = {\bf 0}$).  Finally,
Algorithm~\ref{alg:phaseRetrieval} is fast, with a runtime
complexity that is near-linear in $d$.  We are now ready to begin
proving Theorem~\ref{thm:RecoverArbVecs}.

\subsection{Well Conditioned Measurements}
\label{sec:CondNumBounds}

In this section we develop a set of deterministic measurements $M \in \mathbbm{C}^{D \times d}$ that lead to well conditioned block circulant matrices $M' \in \mathbbm{C}^{D \times D}$ in \eqref{equ:LinProb}.  To begin, we choose $a \in [4, \infty)$ and then set our local correlation masks to be
\begin{equation}
(\m_l)_{i} = \left\{ \begin{array}{ll} \frac{\mathbbm{e}^{-i/a}}{\sqrt[4]{2\delta -1}} \cdot \mathbbm{e}^{\frac{2 \pi \mathbbm{i} \cdot (i - 1) \cdot (l - 1)}{2\delta - 1}} & \textrm{if}~i \leq \delta \\ 0 & \textrm{if}~i > \delta\end{array} \right.
\label{equ:MeasDef}
\end{equation}
for $1 \leq l \leq 2\delta - 1$, and $1 \leq i \leq d$.  This leads to blocks $M'_{l} \in \mathbbm{C}^{(2\delta - 1) \times (2\delta - 1)}$ from \eqref{eq:M'DEF} with entries given by
\begin{equation*}
    \resizebox{\textwidth}{!}{$
(M'_l)_{i,j} := \left\{ \begin{array}{ll} (\m_i)_{l} \overline{(\m_i)}_{j+l-1}  ~=~  \frac{\mathbbm{e}^{-(2l+j-1)/ a}}{\sqrt{2\delta - 1}} \cdot \mathbbm{e}^{-\frac{2 \pi \mathbbm{i} \cdot (i - 1) \cdot (j-1)}{2\delta - 1}} & \textrm{if}~ 1 \leq j \leq \delta - l +1 \\
 0 & \textrm{if}~ \delta - l + 2 \leq j \leq 2\delta - l -1 \\  (\m_i)_{l+1} \overline{(\m_i)}_{l + j - 2 \delta + 1} ~=~  \frac{\mathbbm{e}^{-(2l+j-2(\delta - 1))/a}}{\sqrt{2\delta - 1}} \cdot \mathbbm{e}^{-\frac{2 \pi \mathbbm{i} \cdot (i - 1) \cdot (j - 2 \delta)}{2\delta - 1}}  & \textrm{if}~2 \delta - l \leq j \leq 2\delta - 1, ~l < \delta \\
  0 & \textrm{if}~ j > 1, ~\textrm{and}~ l = \delta  \end{array}
  \right..$}
\end{equation*} 
We will now begin to bound the condition number of this block circulant matrix, $M'$, by block diagonalizing it via \eqref{equ:Diagonalized}.  

Considering the entries of each $J_{k}  \in \mathbbm{C}^{(2\delta-1) \times (2\delta-1)}$ from \eqref{def:JkMat} results in two cases.  First, suppose that $1 \leq j \leq \delta$.  In this case one can see that
\begin{align}
\left( J_{k} \right)_{i,j} ~=&~
\frac{\mathbbm{e}^{(1-j)/a}}{\sqrt{2\delta - 1}} \cdot
\mathbbm{e}^{\frac{-2 \pi \mathbbm{i} \cdot (i - 1) \cdot (j-1)}{2\delta
- 1}}  \cdot \sum_{l = 1}^{\delta-j+1} \mathbbm{e}^{-2l/a} \cdot
\mathbbm{e}^{\frac{2 \pi \mathbbm{i} \cdot (k-1) \cdot (l-1)}{d}}, \\ 
~=&~\frac{\mathbbm{e}^{-(j+1)/a}}{\sqrt{2\delta - 1}} \cdot
\mathbbm{e}^{\frac{-2 \pi \mathbbm{i} \cdot (i - 1) \cdot (j-1)}{2\delta
- 1}}  \cdot \frac{1 - \mathbbm{e}^{-2(\delta-j + 1)/a} \cdot \mathbbm{e}^{\frac{2
  \pi \mathbbm{i} \cdot (k-1) \cdot (\delta-j + 1)}{d}}}{1 -
    \mathbbm{e}^{-2/a} \cdot \mathbbm{e}^{\frac{2 \pi \mathbbm{i} \cdot
  (k-1)}{d}}} .
\end{align}
Second, suppose that $\delta + 1 \leq j \leq 2 \delta - 1$.  In this case one can see that 
\begin{align}
\left( J_{k} \right)_{i,j} ~=&~
\frac{\mathbbm{e}^{-(j-2(\delta-1))/a}}{\sqrt{2\delta - 1}} \cdot
\mathbbm{e}^{\frac{-2 \pi \mathbbm{i} \cdot (i - 1) \cdot
(j-2\delta)}{2\delta - 1}}  \cdot \sum_{l = 2 \delta - j}^{\delta-1}
\mathbbm{e}^{-2l/a} \cdot \mathbbm{e}^{\frac{2 \pi \mathbbm{i} \cdot
(k-1) \cdot (l-1)}{d}}, \\ 
~=&~\frac{\mathbbm{e}^{-(2(\delta + 1) - j)/a}}{\sqrt{2\delta - 1}}
\cdot \mathbbm{e}^{\frac{-2 \pi \mathbbm{i} \cdot (i - 1) \cdot
(j-1)}{2\delta - 1}}  \cdot \mathbbm{e}^{\frac{2 \pi
  \mathbbm{i} \cdot (k-1) (2 \delta - j-1)}{d}} \cdot \frac{1 -
    \mathbbm{e}^{-2(j-\delta)/a} \cdot \mathbbm{e}^{\frac{2 \pi
      \mathbbm{i} \cdot (k-1) \cdot (j - \delta)}{d}}}{1 -
        \mathbbm{e}^{-2/a} \cdot \mathbbm{e}^{\frac{2 \pi \mathbbm{i}
      \cdot (k-1)}{d}}} .
\end{align}

Let $F_\alpha \in \mathbbm{C}^{\alpha \times \alpha}$ be the unitary $\alpha \times \alpha$ discrete Fourier transform matrix.  Defining 
\begin{equation*}
    \resizebox{\textwidth}{!}{$
s_{k,j} ~:= \left\{ \begin{array}{ll}  \displaystyle \mathbbm{e}^{-(j+1)/a} \cdot
  \frac{1 - \mathbbm{e}^{-2(\delta-j + 1)/a} \cdot \mathbbm{e}^{2 \pi
    \mathbbm{i} \cdot (k-1) \cdot (\delta-j + 1)/d}}{1 -
      \mathbbm{e}^{-2/a} \cdot \mathbbm{e}^{2 \pi \mathbbm{i} \cdot
    (k-1)/d}}  & \textrm{if}~ 1 \leq j \leq \delta\\
\displaystyle \mathbbm{e}^{-(2(\delta + 1) - j)/a} \cdot \mathbbm{e}^{2
  \pi \mathbbm{i} \cdot (k-1) (2 \delta - j-1)/d} \cdot \frac{1 -
    \mathbbm{e}^{-2(j-\delta)/a} \cdot \mathbbm{e}^{2 \pi \mathbbm{i}
  \cdot (k-1) \cdot (j - \delta)/d}}{1 - \mathbbm{e}^{-2/a} \cdot
  \mathbbm{e}^{2 \pi \mathbbm{i} \cdot (k-1)/d}} & \textrm{if}~
  \delta + 1 \leq j \leq 2 \delta - 1 \end{array} \right.,$}
\end{equation*} 
we now have that 
\begin{equation}
J_k =  ~F_{2\delta - 1} ~\left( \begin{array}{llll}  
s_{k,1} & 0 & \dots & 0\\
0 & s_{k,2} & 0 & \dots\\
0 & 0 & \ddots & 0 \\
0 & \dots & 0 & s_{k,2\delta - 1}\\
\end{array} \right). 
\label{equ:JExp}
\end{equation}
Note that the condition number of $J$, and therefore of $M'$, will be dictated by the singular values of these $J_k$ matrices.  Thus, we will continue by developing bounds for the singular values of each $J_{k}  \in \mathbbm{C}^{(2\delta-1) \times (2\delta-1)}$.

The fact that $F_{2\delta - 1}$ is unitary implies that
\begin{equation}
\min_{j \in [2\delta-1]} \left| s_{k,j}  \right| \leq  \sigma_{2\delta - 1} \left( J_{k} \right)  \leq  \sigma_1 \left( J_{k} \right) \leq  \max_{j \in [2\delta-1]} \left| s_{k,j}  \right|
\end{equation}
for all $k \in [d]$.  Thus, we will now devote ourselves to bounding the maximum and minimum values of $| s_{k,j} |$ from above and below, respectively, over all $k \in [d]$ and $j \in [2 \delta - 1]$.  These bounds will then collectively yield an upper bound on the condition number of our block circulant measurement matrix $M'$.  The following simple technical lemmas will be useful.

\begin{lem}
Let $x \in [2, \infty)$.  Then, $1 - \mathbbm{e}^{-1 / x} ~>~ \frac{2-\mathbbm{e}^{1 / x}}{x} ~\geq~ \frac{2 - \sqrt{\mathbbm{e}}}{x} > \frac{7}{20 \cdot x}$.  
\label{lem:DenBound}
\end{lem}

\noindent \textit{Proof:}  Note that $1 - \mathbbm{e}^{-1 / x}~=~ \sum^{\infty}_{n=1} \frac{(-1)^{n+1}}{x^n n!}~>~ \frac{1}{x} \cdot \left( 2 - \sum^{\infty}_{n=0} \frac{1}{x^n (n+1)!} \right)  ~>~ \frac{2-\mathbbm{e}^{1 / x}}{x}$.  Furthermore, the numerator is a monotonically increasing function of $x$.  \hfill $\qed$\\

\begin{lem}
Let $a,b,c \in \mathbbm{R}^+$, and $f: \mathbbm{R} \rightarrow \mathbbm{R}$ below.  Then,
\begin{enumerate}
\item $f(x) = b \cdot \mathbbm{e}^{-x / a} \left( 1 + c \cdot \mathbbm{e}^{2x / a} \right)$ has a unique global minimum at $x = -\frac{a}{2} \ln (c)$, and
\item $f(x) = b \cdot \mathbbm{e}^{-x / a} \left( 1 - c \cdot \mathbbm{e}^{2x / a} \right)$ is monotonically decreasing.
\end{enumerate}
\label{lem:BasicCalc}
\end{lem}

\noindent \textit{Proof:}  In either case we have that $f'(x) = -\frac{b}{a}\cdot \mathbbm{e}^{-x / a} \pm \frac{bc}{a} \cdot \mathbbm{e}^{x / a}$, and $f''(x) = \frac{b}{a^2}\cdot \mathbbm{e}^{-x / a} \pm \frac{bc}{a^2} \cdot \mathbbm{e}^{x / a}.$  For $(1)$ we have a single critical point at $x = -\frac{a}{2} \ln (c)$, which is a global minimum since $f''(x) > 0  ~\forall x \in \mathbbm{R}$.  For $(2)$ we have $f'(x) < 0$ for all $x \in \mathbbm{R}$. \hfill $\qed$\\

Note that 
\begin{equation}
    \resizebox{0.925 \textwidth}{!}{$
\left| s_{k,j}  \right| = \left\{ \begin{array}{ll}  
  \displaystyle \mathbbm{e}^{-(j+1)/a} \cdot
  \sqrt{\frac{1+\mathbbm{e}^{-4(\delta - j + 1)/a} - 2
  \mathbbm{e}^{-2(\delta-j+1)/a} \cos\left( 2 \pi \cdot [\delta - j + 1]
\cdot (k-1) / d \right)}{1+\mathbbm{e}^{-4 / a} - 2 \mathbbm{e}^{-2 / a}
\cos\left( 2 \pi (k-1) /d \right)}}  & \textrm{if}~ 1 \leq j \leq \delta\\
\displaystyle \mathbbm{e}^{-(2(\delta + 1) - j)/a} \cdot
\sqrt{\frac{1+\mathbbm{e}^{-4(j - \delta)/a} - 2 \mathbbm{e}^{-2(j -
\delta)/a} \cos\left( 2 \pi \cdot [j - \delta] \cdot (k-1) / d
\right)}{1+\mathbbm{e}^{-4 / a} - 2 \mathbbm{e}^{-2 / a} \cos\left( 2
\pi (k-1) /d \right)}}   & \textrm{if}~ \delta + 1 \leq j \leq 2
\delta - 1 \end{array} \right..$}
\end{equation}
Fix $k \in [d]$.  When $1 \leq j \leq \delta$ we have 
\begin{equation}
\max_{j \in [\delta]} \left| s_{k,j}  \right| \leq \max_{j \in [\delta]} \left( \mathbbm{e}^{-(j+1)/a} \cdot \frac{1+ \mathbbm{e}^{-2(\delta + 1 - j) / a}}{1 - \mathbbm{e}^{-2 / a}} \right) \leq \frac{\mathbbm{e}^{-2 / a}(1+\mathbbm{e}^{-2 \delta / a})}{1 - \mathbbm{e}^{-2 / a}},
\label{equ:UB1}
\end{equation}
where the second inequality follows from part one of Lemma~\ref{lem:BasicCalc}.  When $ \delta + 1 \leq j \leq 2 \delta - 1$ we have 
\begin{equation}
\max_{j \in [2 \delta-1] \setminus [\delta] } \left| s_{k,j}  \right| \leq \max_{j \in [2\delta-1] \setminus [\delta]} \left( \mathbbm{e}^{-(2(\delta + 1) - j)/a} \cdot \frac{1+ \mathbbm{e}^{-2(j - \delta)/a}}{1 - \mathbbm{e}^{-2 / a}} \right) \leq \frac{\mathbbm{e}^{-3 / a}(1+\mathbbm{e}^{-2(\delta-1) / a})}{1 - \mathbbm{e}^{-2 / a}},
\label{equ:UB2}
\end{equation}
where the second inequality again follows from part one of Lemma~\ref{lem:BasicCalc}.  Finally, combining $\eqref{equ:UB1}$ and $\eqref{equ:UB2}$ one can see that
\begin{equation}
\sigma_1 \left( J_{k} \right) \leq \frac{\mathbbm{e}^{-2 / a}(1+\mathbbm{e}^{-2 \delta / a})}{1 - \mathbbm{e}^{-2 / a}} ~<~ a \cdot \frac{\mathbbm{e}^{-2 / a}(1+\mathbbm{e}^{-2 \delta / a})}{2(2 - \mathbbm{e}^{2 / a})} ~<~ a \cdot \frac{20 \mathbbm{e}^{-2 / a}}{7} ~<~ 3 a \cdot \mathbbm{e}^{-2 / a},
\label{equ:UB}
\end{equation}
where the second inequality follows from Lemma~\ref{lem:DenBound} with $a \in [4,\infty)$.

Turning our attention to the lower bound, we note that part two of Lemma~\ref{lem:BasicCalc} implies that
\begin{equation}
\min_{j \in [\delta]} \left| s_{k,j}  \right| \geq \min_{j \in [\delta]} \left( \mathbbm{e}^{-(j+1)/a} \cdot \frac{1- \mathbbm{e}^{-2(\delta + 1 - j) / a}}{1 + \mathbbm{e}^{-2 / a}} \right) \geq \frac{\mathbbm{e}^{-(\delta + 1) / a}(1- \mathbbm{e}^{-2 / a})}{1 + \mathbbm{e}^{-2 / a}}.
\label{equ:LB1}
\end{equation}
Similarly, part two of Lemma~\ref{lem:BasicCalc} also ensures that
\begin{equation}
\min_{j \in [2 \delta - 1] \setminus [\delta]} \left| s_{k,j}  \right| \geq \min_{j \in [2 \delta - 1] \setminus [\delta]} \left( \mathbbm{e}^{-(2(\delta + 1) - j)/a} \cdot \frac{1- \mathbbm{e}^{-2(j - \delta)/a}}{1 + \mathbbm{e}^{-2 / a}} \right) \geq \frac{\mathbbm{e}^{-(\delta + 1) / a}(1- \mathbbm{e}^{-2 / a})}{1 + \mathbbm{e}^{-2 / a}}.
\label{equ:LB2}
\end{equation}
Combining $\eqref{equ:LB1}$ and $\eqref{equ:LB2}$ we see that
\begin{equation}
\sigma_{2 \delta - 1} \left( J_{k} \right) \geq \frac{\mathbbm{e}^{-(\delta + 1) / a}(1- \mathbbm{e}^{-2 / a})}{1 + \mathbbm{e}^{-2 / a}} ~>~ \frac{7}{20 a} \cdot \mathbbm{e}^{-(\delta + 1) / a},
\label{equ:LB}
\end{equation}
where the second inequality follows from Lemma~\ref{lem:DenBound} with $a \in [4,\infty)$.  We are now equipped to prove the main theorem of this section.

\begin{thm}
Define $M' \in \mathbbm{C}^{D \times D}$ via \eqref{equ:MeasDef} with $a:= \max \left\{ 4,~\frac{\delta - 1}{2} \right\}$.  Then, $$\kappa \left( M' \right) ~<~ \max \left\{ 144 \mathbbm{e}^2,~\frac{9 \mathbbm{e}^2}{4} \cdot (\delta - 1)^2 \right \}.$$
\label{thm:WellCondMeas}
\end{thm}

\noindent \textit{Proof:}  We have from \eqref{equ:UB} and \eqref{equ:LB} that
\begin{equation}
\kappa \left( M' \right) = \frac{\sigma_1 \left( M' \right)}{\sigma_{D} \left( M' \right)} ~=~ \frac{\sigma_1 \left( J \right)}{\sigma_{D} \left( J \right)} ~\leq~ \frac{\max_{k \in [d]} \sigma_1 \left( J_k \right)}{\min_{k \in [d]}\sigma_{2 \delta-1} \left( J_k \right)} ~<~ 9 a^2 \cdot \mathbbm{e}^{(\delta - 1) / a}.
\end{equation}
Minimizing the rightmost upper bound as a function of $a$ yields the stated result.  \hfill $\qed$\\

Theorem~\ref{thm:WellCondMeas} guarantees the existence of measurements which allow for the approximation of the phase difference vector ${\bf y} \in \mathbbm{C}^D$ defined in \eqref{equ:Defy}.  In the next three subsections we analyze the approximation of ${\bf x} \in \mathbbm{C}^d$ from ${\bf \tilde{y}} \approx \y$ via the techniques discussed in \S \ref{sec:Alg}.

\subsection{A Recovery Guarantee for Flat Vectors}
\label{sec:RecoverFlatVecs}

As mentioned above, Algorithm~\ref{alg:angSync} implicitly assumes that ${\bf \tilde{y}} \in \mathbbm{C}^D$ does not contain any strings of $\delta - 1$ consecutive entires (mod $d$) all with very small magnitudes.  In this section we will demonstrate that a general class of non-sparse vectors $\x \in \mathbbm{C}^d$ lead to such ${\bf \tilde{y}}$ whenever noise levels are low enough.  As a result, we will prove deterministic approximation guarantees for all sufficiently non-sparse vectors $\x \in \mathbbm{C}^d$ in the high signal-to-noise setting.  More specifically, we will utilize the following concrete characterization of flatness (i.e., non-sparsity) for $\x$ hereafter.

\begin{Def}
Let $m \in [d]$.  A vector $\x \in \mathbbm{C}^d$ will be called $m$-flat if can be partitioned into at least $\left \lfloor \frac{d}{m} \right \rfloor$ blocks of consecutive entries such that:
\begin{enumerate}
\item Every block contains either $m$ or $m+1$ neighboring entries of $\x$, and
\item Every such block of entries contains at least one entry whose magnitude is $\geq \frac{\| \x  \|_2}{2\sqrt{d}}$.
\end{enumerate}
\end{Def}

The following simple lemma proves that Algorithm~\ref{alg:phaseRetrieval} accurately estimates the magnitude of every entry of $\x$.  Both it and the next lemma use the notation from \eqref{eq:InvertMeasDefs} and \eqref{eq:Defk(i,j)}, and implicitly assume the invertibility of $M'$.

\begin{lem}
Let $|\tilde{x}_j| = \sqrt{|\tilde{y}_{k(j,j)}|} = \sqrt{|x_{j}\overline{x}_{j} + \tilde{n}_{k(j,j)}|}$ be the estimate of $|x_j|$ utilized in line 3 of Algorithm~\ref{alg:phaseRetrieval}, where ${\bf \tilde{n}} := (M')^{-1}P \n$ .  Then $\left| |x_j| - |\tilde{x}_j| \right|^2 ~\leq~ 3 \| {\bf \tilde{n}} \|_{\infty}$.
\label{lem:GoodMagEsts}
\end{lem}

\noindent \textit{Proof:}  Let $a:= |x_j| \in \mathbbm{R}^+$ and $\epsilon := \tilde{n}_{k(j,j)} \in \mathbbm{C}$.  We upper bound the righthand side of 
$$\left| |x_j| - |\tilde{x}_j| \right|^2 = a^2 + |a^2 + \epsilon | - 2a \sqrt{|a^2 + \epsilon|} \leq 2a^2 + |\epsilon | - 2a^2 \sqrt{\left|1 + \frac{\epsilon}{a^2}\right|} $$
by considering two cases:  If $a^2 \leq |\epsilon|$ we may bound the rightmost negative term by zero in order to obtain the desired result.  If $a^2 > |\epsilon|$ then
$$-2 a^2 \sqrt{\left| 1 + \frac{\epsilon}{a^2} \right|} ~\leq~ -2 a^2 \sqrt{1 - \frac{|\epsilon|}{a^2}} ~\leq~ -2 a^2 \left( 1 - \frac{|\epsilon|}{a^2} \right)$$
so that
$$\left| |x_j| - |\tilde{x}_j| \right|^2 \leq 2a^2 + |\epsilon | -2 a^2 \left( 1 - \frac{|\epsilon|}{a^2} \right) = 3 |\epsilon|$$
again holds as desired.\hfill $\qed$\\

The next lemma proves that Algorithm~\ref{alg:angSync} will accurately estimate the phases of all entires $x_j$ whose magnitudes are sufficiently large relative to the noise level.  This lemma requires that the vector $\x$ be $m$-flat for a sufficiently small $m$.

\begin{lem}
Suppose that $\x \in \mathbbm{C}^d$ is $\left \lfloor \frac{\delta - 3}{2} \right \rfloor$-flat with $d > 2$ and $\| \x \|^2_2 \geq 26 ~ d^2 ~\| {\bf \tilde{n}} \|_\infty$, where ${\bf \tilde{n}} = (M')^{-1}P \n$ as in \eqref{eq:InvertMeasDefs}.  Let $k \in [d]$ be such that $|x_k|^2 \geq \frac{5}{2} ~d~ \|{\bf \tilde{n}} \|_\infty$, and $\phi_a \in [0, 2 \pi]$ be the true phase of the entry $x_a$ of $\x$ chosen in line 2 of Algorithm~\ref{alg:angSync}.  Then, Algorithm~\ref{alg:angSync} will produce an estimate $\tilde{\phi}_k$ of the true phase of $x_k$, $\phi_k := \arg(x_k)$, satisfying 

$$\left| \mathbbm{e}^{\mathbbm{i} \tilde{\phi}_k} - \mathbbm{e}^{\mathbbm{i} ( \phi_k - \phi_a ) } \right|^2 \leq \frac{10 ~d^2~ \| {\bf \tilde{n}} \|_\infty}{\| \x \|^2_2}.$$
\label{lem:Alg2Error}
\end{lem}

\noindent \textit{Proof:} Let $\alpha := \frac{\| \x \|_2}{d \sqrt{10  \| {\bf \tilde{n}} \|_\infty}}$ and $\beta := \frac{\| \x  \|_2}{2\sqrt{d}}$.  Note that $\alpha > 1$, $|x_k| \geq \frac{\beta}{\alpha}$, $\| {\bf \tilde{n}} \|_{\infty} < \frac{1}{4} \frac{\beta^2}{d \alpha}$, and $| x_a | > \frac{\sqrt{3}}{2} \beta$ all hold.  Furthermore, the fact that $\x$ is $\left \lfloor \frac{\delta - 3}{2} \right \rfloor$-flat guarantees the existence of a sequence of entries $x_{b_1}, \dots, x_{b_q}$ such that: 
\begin{enumerate}
\item $0 < |(b_{l+1} - b_l) \textrm{~mod~d}| \leq \delta - 1$ for all $l \in [q-1]$, 
\item $0 \leq \max \{ |(b_1 - a) \textrm{~mod~d}|, |(k-b_q)  \textrm{~mod~d}| \} \leq \delta - 1$, and
\item $|x_{b_l}| \geq \beta$ for all $l \in [q]$.
\end{enumerate}
Thus, line 18 of Algorithm~\ref{alg:angSync} is guaranteed to identify a sequence of entries $x_{j_1}, \dots, x_{j_p}$ such that:
\begin{enumerate}
\item $0 < |(j_{l+1} - j_l) \textrm{~mod~d}| \leq \delta - 1$ for all $l \in [q-1]$, 
\item $0 \leq \max \{ |(j_1 - a) \textrm{~mod~d}|, |(k-j_p)  \textrm{~mod~d}| \} \leq \delta - 1$, and
\item $|x_{j_l}| > \frac{\sqrt{3}}{2} \beta$ for all $l \in [q]$.
\end{enumerate}
These entries contribute to the estimate $\tilde{\phi}_k$ of $\phi_{k} - \phi_{a}$ given by Algorithm~\ref{alg:angSync},
\begin{equation}
\tilde{\phi}_k ~=~ \tilde{\phi}_{k,j_p} + \left( \sum^{p-1}_{l = 1} \tilde{\phi}_{j_{l+1}, j_l} \right) + \tilde{\phi}_{j_{1}, a} ~\approx~(\phi_{k} - \phi_{j_{p}}) + \left( \sum^{p-1}_{l = 1} \phi_{j_{l+1}} - \phi_{j_{l}} \right) + (\phi_{j_{1}} - \phi_{a}) = \phi_{k} - \phi_{a},
\label{eq:Alg2Error}
\end{equation}
where $\tilde{\phi}_{j_{l+1}, j_l} := \frac{1}{2} \left( \arg \left( \tilde{y}_{k(j_{l+1},j_l)} \right) - \arg \left( \tilde{y}_{k(j_l,j_{l+1})} \right) \right)$ is used as an estimate in line 13 of $\phi_{j_{l+1},j_l} := \arg\left( y_{k(j_{l+1},j_l)} \right) = \phi_{j_{l+1}} - \phi_{j_{l}}$.

To bound the approximation error in \eqref{eq:Alg2Error} we note that the law of sines implies that
$$\sin \left( | \tilde{\phi}_{j_{l+1}, j_l} - \phi_{j_{l+1},j_l} | \right) \leq \left| \sin \left( \tilde{\phi}_{j_{l+1}, j_l} - \phi_{j_{l+1},j_l} \right) \right| \leq \frac{\| {\bf \tilde{n}} \|_{\infty}}{|\tilde{y}_{l'}|},$$
where $\tilde{y}_{l'} = x_{j_{l+1}} \overline{x}_{j_l} + \tilde{n}_{l'}$ is, without loss of generality, the smaller of the two measurements (in magnitude) that contribute to the estimate $\tilde{\phi}_{j_{l+1}, j_l}$.  Using the facts from the previous paragraph concerning $|x_a|$, $|x_k|$, and $|x_{j_l}|$ for all $l \in [p]$, one can show that every measurement $\tilde{y}_{l'}$ that contributes to an estimate $\tilde{\phi}_{j_{l+1}, j_l}$ used in \eqref{eq:Alg2Error} will have $|\tilde{y}_{l'}| > \frac{\beta^2}{2 \alpha}$.  Hence, we have that
$$\sin \left( | \tilde{\phi}_{j_{l+1}, j_l} - \phi_{j_{l+1},j_l} | \right) \leq \frac{2\alpha}{\beta^2}\| {\bf \tilde{n}} \|_{\infty} < \frac{1}{4}$$
holds for all $l \in [p]$.  As a consequence, we can infer that
$$\left| \tilde{\phi}_{j_{l+1}, j_l} - (\phi_{j_{l+1}} - \phi_{j_{l}} ) \right| = \left| \tilde{\phi}_{j_{l+1}, j_l} - \phi_{j_{l+1},j_l} \right| \leq \frac{5}{4} \sin \left( | \tilde{\phi}_{j_{l+1}, j_l} - \phi_{j_{l+1},j_l} | \right) \leq \frac{5}{2} \frac{\alpha}{\beta^2}\| {\bf \tilde{n}} \|_{\infty}$$
also holds for all all $l \in [p]$.  Combining this with \eqref{eq:Alg2Error} we learn that 
$$\left| \tilde{\phi}_k - \phi_k + \phi_a \right| \leq d \frac{5}{2} \frac{\alpha}{\beta^2}\| {\bf \tilde{n}} \|_{\infty} < \frac{5}{8}.$$
This, in turn, implies that
$$\left| \mathbbm{e}^{\mathbbm{i} \tilde{\phi}_k} - \mathbbm{e}^{\mathbbm{i} ( \phi_k - \phi_a ) } \right| = \left| \mathbbm{e}^{\mathbbm{i} (\tilde{\phi}_k - \phi_k + \phi_a)} - 1 \right| = 2 \sin \left( \frac{1}{2} \left| \tilde{\phi}_k - \phi_k + \phi_a \right| \right) \leq  \left| \tilde{\phi}_k - \phi_k + \phi_a \right| \leq \frac{5d \alpha}{2 \beta^2} \| {\bf \tilde{n}} \|_{\infty}$$
proving the lemma. \hfill $\qed$\\

We are now prepared to prove the main result of this section.

\begin{thm}
There exist fixed universal constants $C, C' \in \mathbbm{R}^+$ such that following holds:  Let $M \in \mathbbm{C}^{D \times d}$ be defined as in \S\ref{sec:CondNumBounds}, and suppose that $\x \in \mathbbm{C}^d$ is $\left \lfloor \frac{\delta - 3}{2} \right \rfloor$-flat with $d > 2$ and $\| \x \|^2_2 \geq C ~ (\delta - 1) d^2 ~\| {\bf n} \|_2$.  Then, Algorithm~\ref{alg:phaseRetrieval} is guaranteed to recover an ${\bf \tilde{x}} \in \mathbbm{C}^d$ with 
\begin{equation}
\min_{\theta \in [0, 2 \pi]} \left\| {\bf x}  - \mathbbm{e}^{\mathbbm{i} \theta} {\bf \tilde{x}} \right\|^2_2 ~\leq~ C' d^2 (\delta - 1) \| {\bf n} \|_2
\label{eq:Recovery!}
\end{equation}
when given arbitrarily noisy input measurements $\b = \left \vert M \x \right \vert^2 + \n \in \mathbbm{R}^D$ as per \eqref{eq:PR_problem}.  Furthermore, Algorithm~\ref{alg:phaseRetrieval} requires just $\mathcal{O}(\delta \cdot d \log d)$ operations for this choice of $M \in \mathbbm{C}^{D \times d}$.
\label{thm:RecoverFlatVecs}
\end{thm}

\noindent \textit{Proof:}  
First, note that $M'$ will be invertible by Theorem~\ref{thm:WellCondMeas}.  Thus, we may set ${\bf \tilde{n}} = (M')^{-1}P \n$.  Furthermore, the proof of Theorem~\ref{thm:WellCondMeas} tells us that there exists an explicit universal constant $C'' \in \mathbbm{R}^+$ such that
\begin{equation}
C'' (\delta - 1) \| {\bf n} \|_2 \geq \frac{\| {\bf n} \|_2}{\sigma_{D} \left( M' \right)} \geq  \| {\bf \tilde{n}} \|_2 \geq  \| {\bf \tilde{n}} \|_\infty.
\label{equ:CondBound4Thm} 
\end{equation}
Setting $C = 26 C''$ now allows us to verify that
$$\| \x \|^2_2 \geq 26 C'' ~ (\delta - 1) d^2 ~\| {\bf n} \|_2 \geq 26 ~d^2~ \| {\bf \tilde{n}} \|_\infty.$$
Hence, we will be able to apply both Lemmas~\ref{lem:GoodMagEsts} and~\ref{lem:Alg2Error} as desired.

Let $\phi_a \in [0, 2 \pi]$ be as in Lemma~\ref{lem:Alg2Error}, and let ${\bf p}, {\bf \tilde{p}} \in \mathbbm{C}^d$ be the vectors of phases of the entries of $\mathbbm{e}^{-\mathbbm{i} \phi_a }\x$ and ${\bf \tilde{x}}$, respectively, so that $p_j = \mathbbm{e}^{\mathbbm{i} ( \phi_j - \phi_a ) }$ and $\tilde{p}_j = \mathbbm{e}^{\mathbbm{i} \tilde{\phi}_j }$ hold for all $j \in [d]$.  Similarly, recall that $|\x|, |{\bf \tilde{x}}| \in \mathbbm{R}^d$ denote the vectors of magnitudes of the entries of $\x$ and ${\bf \tilde{x}}$, respectively, so that $|x|_j = |x_j|$ and $|\tilde{x}|_j = \sqrt{|\tilde{y}_{k(j,j)}|} = \sqrt{|x_{j}\overline{x}_{j} + \tilde{n}_{k(j,j)}|}$ hold for all $j \in [d]$ (here we are again using the notation from \eqref{eq:InvertMeasDefs} and \eqref{eq:Defk(i,j)}).  Thus, $\mathbbm{e}^{-\mathbbm{i} \phi_a }\x = |\x| \circ {\bf p}$ and ${\bf \tilde{x}} = |{\bf \tilde{x}}| \circ {\bf \tilde{p}}$ both hold, where $\circ$ denotes the entrywise (Hadamard) product.

To obtain \eqref{eq:Recovery!} we bound 
\begin{align}
\| \mathbbm{e}^{-\mathbbm{i} \phi_a }\x - {\bf \tilde{x}} \|_2 ~&\leq~ \left\| |\x| \circ {\bf p} - |\x| \circ {\bf \tilde{p}} \right\|_2 + \left \| |\x| \circ {\bf \tilde{p}}  - |{\bf \tilde{x}}| \circ {\bf \tilde{p}} \right \|_2 \nonumber\\
&=~ \sqrt{\sum^d_{j=1}  \left| x_j \right|^2 \left| \mathbbm{e}^{\mathbbm{i} \tilde{\phi}_j} - \mathbbm{e}^{\mathbbm{i} ( \phi_j - \phi_a ) } \right|^2} + \sqrt{\sum^d_{j=1}  \left| |x_j| - |\tilde{x}_j| \right|^2 } \nonumber \\
&\leq~\sqrt{\sum_{j \in S}  \left| x_j \right|^2 \left| \mathbbm{e}^{\mathbbm{i} \tilde{\phi}_j} - \mathbbm{e}^{\mathbbm{i} ( \phi_j - \phi_a ) } \right|^2} +\sqrt{\sum_{j\in[d]\setminus S}  \left| x_j \right|^2 \left| \mathbbm{e}^{\mathbbm{i} \tilde{\phi}_j} - \mathbbm{e}^{\mathbbm{i} ( \phi_j - \phi_a ) } \right|^2} + \sqrt{ 3d \| {\bf \tilde{n}} \|_{\infty} },
\label{eq:FlatRecPr1}
\end{align}
where $S := \left \{ k \in [d] ~\big |~ |x_k|^2 \geq \frac{5}{2} ~d~ \|{\bf \tilde{n}} \|_\infty \right \}$, and the last inequality results from an application of Lemma~\ref{lem:GoodMagEsts}.  Bounding the first two terms of \eqref{eq:FlatRecPr1} using Lemma~\ref{lem:Alg2Error} and the properties of $S$ we can see that
\begin{align*}
\| \mathbbm{e}^{-\mathbbm{i} \phi_a }\x - {\bf \tilde{x}} \|_2 ~&\leq~  \sqrt{\sum_{j \in S}  \left| x_j \right|^2 \left( \frac{10 ~d^2~ \| {\bf \tilde{n}} \|_\infty}{\| \x \|^2_2} \right)}  ~+~\sqrt{ 10 d^2 \| {\bf \tilde{n}} \|_{\infty} } ~+~ \sqrt{ 3d \| {\bf \tilde{n}} \|_{\infty} }\\
&\leq~  2d\sqrt{10 \| {\bf \tilde{n}} \|_\infty} ~+~ \sqrt{ 3d \| {\bf \tilde{n}} \|_{\infty} } ~\leq~ 3d\sqrt{10 \| {\bf \tilde{n}} \|_\infty}.
\end{align*}
Using \eqref{equ:CondBound4Thm} now allows us to establish that
$$\| \mathbbm{e}^{-\mathbbm{i} \phi_a }\x - {\bf \tilde{x}} \|_2 ~\leq~ 3d\sqrt{10 \| {\bf \tilde{n}} \|_\infty} ~\leq~ 3d\sqrt{10C'' (\delta - 1) \| {\bf n} \|_2}$$
which implies \eqref{eq:Recovery!}.

We finish by noting that the runtime complexity of Algorithm~\ref{alg:phaseRetrieval} simplifies to $\mathcal{O}(\delta \cdot d \log d)$ operations when using the measurements defined in \S\ref{sec:CondNumBounds} because the matrix $J$ also has a simple block-diagonal factorization in this case (recall \eqref{equ:Diagonalized} and \eqref{equ:JExp} in light of \S\ref{sec:AlgRuntimeAnal}).
 \hfill $\qed$\\ 

The following corollary easily follows from Theorem~\ref{thm:RecoverFlatVecs}.

\begin{cor}
Let $M \in \mathbbm{C}^{D \times d}$ be defined as in \S\ref{sec:CondNumBounds}, and suppose that ${\bf x} \in \mathbbm{C}^d$ is $m$-flat for some $m \leq \left \lfloor \frac{\delta - 3}{2} \right \rfloor$.  Then, Algorithm~\ref{alg:phaseRetrieval} will recover an ${\bf \tilde{x}} \in \mathbbm{C}^d$ with ${\bf \tilde{x}} = \mathbbm{e}^{\mathbbm{i} \theta} {\bf x}$ for some $\theta \in [0, 2 \pi]$ when given noiseless input measurements $\b = |  M \x |^2 \in \mathbbm{R}^D$.  Furthermore, Algorithm~\ref{alg:phaseRetrieval} requires just $\mathcal{O}(\delta \cdot d \log d)$ operations in this case.
\end{cor}

Of course, not all vectors are $m$-flat for a suitably small value of $m$.  We will generalize our results to arbitrary vectors $\x$ in the next section.  This will be accomplished by showing that a well chosen random unitary matrix $W$ will have the property that $W \x$ is $m$-flat with high probability.

\subsection{Flattening Arbitrary Vectors with High Probability}
\label{sec:FlattenVecs}

Let $W \in \mathbbm{C}^{d \times d}$ be the random unitary matrix
\begin{equation}
W := P F B,
\label{def:Wdef}
\end{equation}
where $P \in \{ 0,1 \}^{d \times d}$ is a permutation matrix selected uniformly at random from the set of all $d \times d$ permutation matrices, $F$ is the unitary $d \times d$ discrete Fourier transform matrix, and $B \in \left\{ -1,0,1\right\}^{d \times d}$ is a random diagonal matrix with i.i.d. symmetric Bernoulli entries on its diagonal.  For any given $m \in [d]$, one can naturally partition $W$ into $\left \lfloor  \frac{d}{m} \right \rfloor$ blocks of contiguous rows, each of cardinality either $m$ or $m+1$.  This defines the $\left \lfloor  \frac{d}{m} \right \rfloor$ sub-matrices of $W$, $W_1, \dots, W_{d - m\left \lfloor  \frac{d}{m} \right \rfloor} \in \mathbbm{C}^{(m+1) \times d}$ and $W_{d - m\left \lfloor  \frac{d}{m} \right \rfloor + 1}, \dots, W_{\left \lfloor  \frac{d}{m} \right \rfloor} \in \mathbbm{C}^{m \times d}$, by
\begin{equation}
W = \left( \begin{array}{l}  W_1\\ W_2\\ \vdots \\ W_{\left \lfloor  \frac{d}{m} \right \rfloor}\end{array} \right).
\label{equ:subW}
\end{equation}

Note that each renormalized sub-matrix of $W$, $\sqrt{\frac{d}{m}} \cdot W_j$ for $j \in \left[ \left \lfloor  \frac{d}{m} \right \rfloor \right]$, is almost a random sampling matrix \eqref{equ:RandSamp} times a random diagonal Bernoulli matrix.  As a result, Theorems~\ref{thm:KW} and~\ref{thm:FourierBON} suggest that each $\sqrt{\frac{d}{m}} \cdot W_j$ should behave like a JL($m$,$d$,$\epsilon$)-embedding of our signal $\x$ into $\mathbbm{C}^m$ (or $\mathbbm{C}^{m+1}$).  If true, it would then be reasonable to expect that each block of $m$ consecutive entries of $W \x$ should have roughly the same $\ell_2$-norm as one another.  This, in turn, suggests that the random unitary matrix $W$ should effectively flatten $\x$ with high probability, especially when $m$ is small.

Of course, there are several small difficulties that must be addressed before the argument above can be made rigorous.  First, the rows of $F$ contributing to $\sqrt{\frac{d}{m}} \cdot W_j$ are effectively independently sampled uniformly \textit{without replacement} from the set of all rows of $F$ by our choice of $P$.\footnote{Note that $W$ being unitary helps us to be able to guarantee both exact recovery of $\x$ in the noiseless setting, and well behaved approximation of $\x$ in the noisy setting.  If we had chosen rows of $F$ with replacement instead of without replacement in \eqref{def:Wdef} we would not have a unitary (or even invertible) square matrix $W$ with probability $\rightarrow 1/\mathbbm{e}$ as $d \rightarrow \infty$.  If one decides to make $W$ rectangular instead of square simply in order to allow sampling from the rows of $F$ with replacement, then many other small difficulties and inefficiencies result.  Thus, we let $P$ be a randomly selected permutation matrix instead of creating it by putting a one in each row independently at random.}  This means that Theorem~\ref{thm:FourierBON} does not strictly apply in our situation since we can not select any row of $F$ more than once.  Secondly, some care must be taken in order to select the smallest value of $m$ possible in \eqref{equ:subW}, since $W \x$ will become flatter, or less sparse, as $m$ decreases.  As a result, $m$ will effectively provide a theoretical lower bound on the size of $\delta$ that one can utilize and still be guaranteed to accurately recover $W \x$ via our \S \ref{sec:Alg} techniques (recall also \S \ref{sec:RecoverFlatVecs} above).  We are now ready to begin proving our main result concerning $W$.

The following simple lemma will be used in order to help adapt Theorem~\ref{thm:FourierBON} to the situation where the rows of $F$ are sampled uniformly without replacement.

\begin{lem}
Let $m \in \mathbbm{N}$ with $m \leq \sqrt{d}$.  Independently draw $x_1, \dots, x_m$ from $[d]$ uniformly at random with replacement.  Then, $\mathbbm{P} \left[ \left| \{ x_1, \dots, x_m \} \right| = m \right] \geq 1/2$.
\label{lem:BDPsimp}
\end{lem}

\noindent \textit{Proof:}
A short induction argument establishes that 
\begin{equation}
\mathbbm{P} \left[ \left| \{ x_1, \dots, x_m \} \right| = m \right] = \prod^{m-1}_{j=1} \left( 1 - \frac{j}{d} \right) \geq 1 - \sum^{m-1}_{j=1} \frac{j}{d} = 1-\frac{m^2-m}{2d}.
\end{equation}
The result now follows easily via algebraic manipulation.\hfill $\qed$\\

The following corollary of Theorem~\ref{thm:FourierBON} now demonstrates that a random sampling matrix $R'$ formed by sampling a subset of rows of size $m$ uniformly at random from $F$ will still be  RIP($2s$,$\epsilon / C_1$) with high probability.

\begin{cor}
Let $p \in (0,1)$.  Form a random sampling matrix $R' \in \mathbbm{C}^{m \times d}$ by independently sampling $m$ rows from $F$ uniformly without replacement.  If the number of rows, $m$, satisfies both 
\begin{equation}
\sqrt{d} \geq m \geq C_3 \cdot \frac{s \ln^2(8s) \ln(8d) \ln(9 m)}{\epsilon^2}
\end{equation}
and
\begin{equation}
\sqrt{d} \geq m \geq C_4 \cdot \frac{s \log(2/p)}{\epsilon^2},
\end{equation}
then $R'$ will be RIP($2s$,$\epsilon / C_1$) with probability at least $1-p$.
\label{cor:FourierBON}
\end{cor}

\noindent \textit{Proof:}
Let $\mathcal{S} := \{ x_1, \dots, x_m \}$, where each $x_j \in [d]$ is selected independently and uniformly at random from $[d]$ (with replacement).  Similarly, let $\mathcal{S}' \subset [d]$ be a subset of $[d]$ chosen uniformly at random from all subsets of $[d]$ with cardinality $m$ (i.e., let $\mathcal{S}'$ contain $m$ elements sampled independently and uniformly from $[d]$ without replacement).  Furthermore, let $E$ denote the event that the random sampling matrix whose rows from $F$ are $x_1, \dots, x_m$ is not RIP($2s$,$\epsilon / C_1$).  Finally, let $E'$ denote the event that the random sampling matrix whose rows from $F$ are the elements of $\mathcal{S}'$ is not RIP($2s$,$\epsilon / C_1$).  Applying Lemma~\ref{lem:BDPsimp} we can now see that
\begin{equation}
\mathbbm{P} \left[ E \right] ~\geq~ \mathbbm{P} \left[ E ~\big|~ |\mathcal{S}| = m \right] \cdot \mathbbm{P} \left[ |\mathcal{S}| = m \right] ~=~ \mathbbm{P} \left[ E' \right] \cdot \mathbbm{P} \left[ |\mathcal{S}| = m \right] ~\geq~ \frac{1}{2} \cdot \mathbbm{P} \left[ E' \right] .
\end{equation}
The stated result now follows from Theorem~\ref{thm:FourierBON}.\hfill $\qed$\\

We are now ready to prove that $W$ will flatten the signal $\x \in \mathbbm{C}^d$ with high probability provided that $m$ can be chosen appropriately.  We have the following theorem:

\begin{thm}
Let $W \in \mathbbm{C}^{d \times d}$ be formed as per \eqref{def:Wdef} for $d \geq 8$.  Then, $W \x \in \mathbbm{C}^d$ will be $m$-flat with probability at least $1-\frac{1}{m}$ provided that $\sqrt{d} \geq m + 1 \geq C_5 \cdot \ln^2(d) \cdot \ln^3 \left( \ln d \right)$.\footnote{Here $C_5 \in \mathbbm{R}^+$  is a fixed absolute constant.}
\label{thm:flattenx}
\end{thm}

\noindent \textit{Proof:}  Our first goal will be to show that each $W_1, \dots, W_{\left \lfloor  \frac{d}{m} \right \rfloor}$ from \eqref{equ:subW} is a is a rescaled JL($m$,$d$,$1/2$)-embedding of $\{ \x \}$ into $\mathbbm{C}^m$ (or $\mathbbm{C}^{m+1}$).  This will guarantee that each consecutive block of $m$ (or $m+1$) entries of $W \x$ has roughly the same $\ell_2$-norm.  

To achieve this goal we will apply Theorem~\ref{thm:KW} to each $\sqrt{\frac{d}{m}} \cdot W_1, \dots, \sqrt{\frac{d}{m}} \cdot W_{\left \lfloor  \frac{d}{m} \right \rfloor}$ in order to show that each one embeds $\{ \x \}$ into $\mathbbm{C}^m$ (or $\mathbbm{C}^{m+1}$) with probability at least $1 - \frac{1}{2d}$.  The union bound will then imply that $\{ \x \}$ is embedded by all the $\sqrt{\frac{d}{m}} \cdot W_j$ with probability at least $1 - \frac{1}{2m}$.  This argument will go through as long as each $\sqrt{\frac{d}{m}} \cdot W_1 B^{-1}, \dots, \sqrt{\frac{d}{m}} \cdot W_{\left \lfloor  \frac{d}{m} \right \rfloor} B^{-1}$ is RIP($2s$,$1 / 2C_1$) for some $s \geq C_2 \cdot \ln(8d)$.  Hence, we will now focus on determining the range of $m$ which guarantees that all $\left \lfloor  \frac{d}{m} \right \rfloor$ of these matrices are RIP($\lceil 2 C_2 \cdot \ln(8d) \rceil$,$1 / 2C_1$).

To demonstrate that each $\sqrt{\frac{d}{m}} \cdot W_j B^{-1}$ is RIP($\lceil 2 C_2 \cdot \ln(8d) \rceil$,$1 / 2C_1$) with probability at least $1-\frac{1}{2d}$ one may apply Corollary~\ref{cor:FourierBON} with $m$ (or $m+1$) chosen as above (assuming $d \geq 8$).  Another application of the union bound then establishes that all of $\sqrt{\frac{d}{m}} \cdot W_1 B^{-1}, \dots, \sqrt{\frac{d}{m}} \cdot W_{\left \lfloor  \frac{d}{m} \right \rfloor} B^{-1}$ will be RIP($\lceil 2 C_2 \cdot \ln(8d) \rceil$,$1 / 2C_1$) with probability at least $1 - \frac{1}{2m}$.  One final application of the union bound then establishes our first goal:  All of $\sqrt{\frac{d}{m}} \cdot W_1, \dots, \sqrt{\frac{d}{m}} \cdot W_{\left \lfloor  \frac{d}{m} \right \rfloor}$ will be JL($m$,$d$,$1/2$)-embeddings of $\{ \x \}$ with probability at least $1 - \frac{1}{m}$.  

To finish the proof, we now note that $W \x$ will be $m$-flat whenever all $\left \lfloor  \frac{d}{m} \right \rfloor$ of the $\sqrt{\frac{d}{m}} \cdot W_j$ matrices are JL($m$,$d$,$1/2$)-embeddings of $\{ \x \}$.  To see why, suppose that 
$$\frac{1}{2} \| ~\x~ \|_2^2 \leq \frac{d}{m} \| W_j \x \|_2^2 \leq \frac{3}{2} \| ~\x~ \|_2^2.$$
This implies that $\frac{3m}{2d} \| ~\x~ \|_2^2 \geq \| W_j \x \|^2_2 \geq \frac{m}{2d} \| \x \|^2_2$, which can only happen if both of the following hold:  $(i)$ at least one entry of $W_j \x$ has magnitude at least $\frac{ \| \x \|_2}{2\sqrt{d}} = \frac{ \| W\x \|_2}{2\sqrt{d}}$, and $(ii)$ all entires of $W_j \x$ have magnitude less than $\sqrt{\frac{3 m+ 3}{2 d }} \| \x \|_2 = \sqrt{\frac{3 m + 3}{2 d }} \| W\x \|_2$.  This proves the theorem. \hfill $\qed$\\

Theorem~\ref{thm:flattenx} now allows us to alter our measurements so that we can recover arbitrary vectors.  We are now ready to prove Theorem~\ref{thm:RecoverArbVecs}.

\subsection{Proof of Theorem~\ref{thm:RecoverArbVecs}}  We set our measurement matrix $\tilde{M} \in \mathbbm{C}^{D \times d}$ to be $\tilde{M} := MW$ where $M \in \mathbbm{C}^{D \times d}$ is defined as in \S\ref{sec:CondNumBounds}, and $W \in \mathbbm{C}^{d \times d}$ is as defined as in \eqref{def:Wdef}.  Theorem~\ref{thm:flattenx} guarantees that $W \x$ will be $m = \mathcal{O}( \ln^2(d) \cdot \ln^3 \left( \ln d \right))$-flat with probability at least $1-\frac{1}{C_5 \cdot \ln^2(d) \cdot \ln^3 \left( \ln d \right)}$ provided that $d$ is sufficiently large.  Furthermore, if $W \x$ is $m$-flat and $\delta \geq 2m+3$, then Theorem~\ref{thm:RecoverFlatVecs} guarantees that Algorithm~\ref{alg:phaseRetrieval} will recover an ${\bf x'} \in \mathbbm{C}^d$ satisfying  
\begin{equation}
\min_{\theta \in [0, 2 \pi]} \left\| W{\bf x}  - \mathbbm{e}^{\mathbbm{i} \theta} {\bf x'} \right\|_2 \leq C' d^2 (\delta - 1) \| {\bf n} \|_2
\label{eq:MainThmpfHalfway}
\end{equation}
when given the noisy input measurements $\b = | M W\x |^2  + \n \in \mathbbm{R}^D$.  Hence, choosing $\delta = \mathcal{O}( \ln^2(d) \cdot \ln^3 \left( \ln d \right))$ allows us to recover ${\bf x'} \approx W \left(\mathbbm{e}^{\mathbbm{i} \phi} \x \right)$, for some unknown phase $\phi \in [0, 2 \pi],$ with probability at least $1-\frac{1}{C_5 \cdot \ln^2(d) \cdot \ln^3 \left( \ln d \right)}$.\footnote{The probability estimate in Theorem~\ref{thm:RecoverArbVecs} follows immediately with $C = C_5$.}  Setting ${\bf \tilde{x}} = W^* {\bf x'}$ we can see that
\begin{align*}
\min_{\theta \in [0, 2 \pi]} \left\| {\bf x}  - \mathbbm{e}^{\mathbbm{i} \theta} {\bf \tilde{x}} \right\|_2 ~&=~ \min_{\theta \in [0, 2 \pi]} \left\| W^* \left( W{\bf x}  - \mathbbm{e}^{\mathbbm{i} \theta} {\bf x'} \right) \right\|_2 \\
&=~\min_{\theta \in [0, 2 \pi]} \left\| W{\bf x}  - \mathbbm{e}^{\mathbbm{i} \theta} {\bf x'} \right\|_2 ~\leq~ C' d^2 (\delta - 1) \| {\bf n} \|_2
\end{align*}
by \eqref{eq:MainThmpfHalfway} since $W$ is always unitary.

Considering the runtime complexity, we note that ${\bf x'}$ can be obtained in $\mathcal{O}(\delta \cdot d \log d)$ = $\mathcal{O}(d \cdot \ln^3(d) \cdot \ln^3 \left( \ln d \right))$ operations by Theorem~\ref{thm:RecoverFlatVecs}.  Computing $W^* {\bf x'}$ can then be done in $\mathcal{O}(d \log d)$ operations via an inverse fast Fourier transform.  The stated runtime complexity follows.\\

It is interesting to note that alternate constructions of flattening matrices, $W$, with fast inverse-matrix vector multiplies can also be created by using sparse Johnson-Lindenstrauss embedding matrices in the place of our Fourier-based matrices (see, e.g., \cite{bourgain2013toward}).  Thus, one has several choices of matrices $W$ to use in concert with a given block-circulant measurement matrix $M$ in principle.

\section{Empirical Evaluation}
\label{sec:Eval}
We now present numerical results demonstrating the efficiency and
robustness of the {\em BlockPR} algorithm. We test our algorithm
on i.i.d. zero-mean complex random Gaussian test signals. To test
noise robustness, we add i.i.d random Gaussian noise to the
squared magnitude measurements at desired signal to noise ratios
(SNRs).  In particular, the noise vector $\n \in \mathbbm R^D$ in
(\ref{eq:PR_problem}) is chosen to be i.i.d. $\mathcal N(\mathbf
0,\sigma^2\mathcal I_D)$. The variance $\sigma^2$ is chosen such
that 
\[ \mbox{SNR (dB)} = 10 \log_{10} \left( \frac{\| M\x\|_2^2}
		{D \, \sigma^2}\right). \]
Errors in the recovered signal are also reported in dB with 
\[ \mbox{Error (dB)} = 10 \log_{10} \left( \frac{\|\tilde \x - 
	\x \|_2^2}{\| \x \|^2_2}\right), \]
where $\tilde \x$ denotes the recovered signal. Matlab code used
to generate the numerical results is freely available at
\cite{bitbucket_BlockPR}. 
\begin{figure}[hbtp]
\centering
\begin{subfigure}[b]{0.475\textwidth}
\includegraphics[clip=true, trim = 0.5in 2in 0.25in 2.45in,scale=0.45]{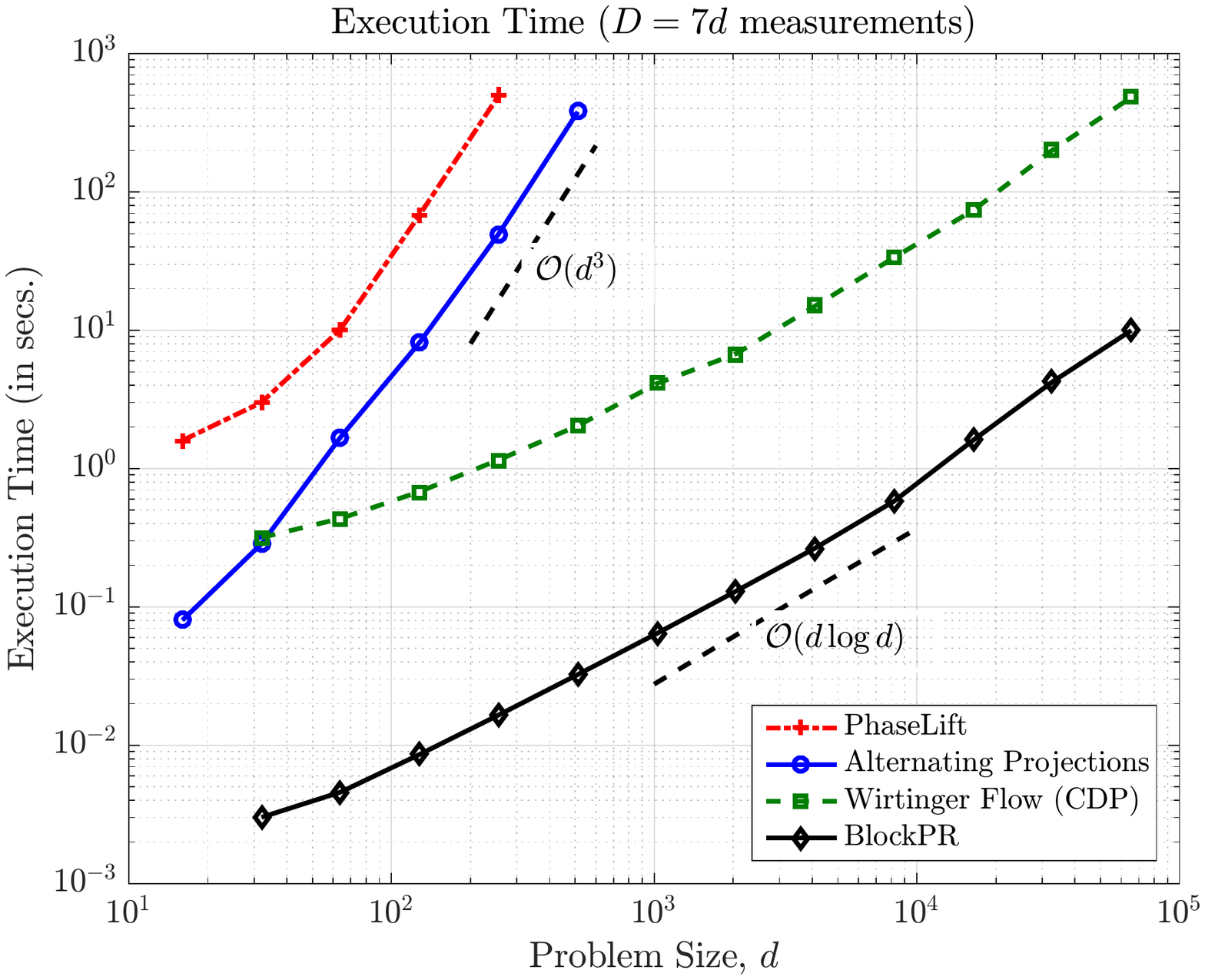}
\caption{Execution time -- Phase Retrieval from $D=7d$ 
measurements.}
\label{fig:highsnr_runtime}
\end{subfigure}
\hfill
\begin{subfigure}[b]{0.475\textwidth}
\includegraphics[clip=true, trim = 0.5in 2in 0.25in 2.45in,
scale=0.45]{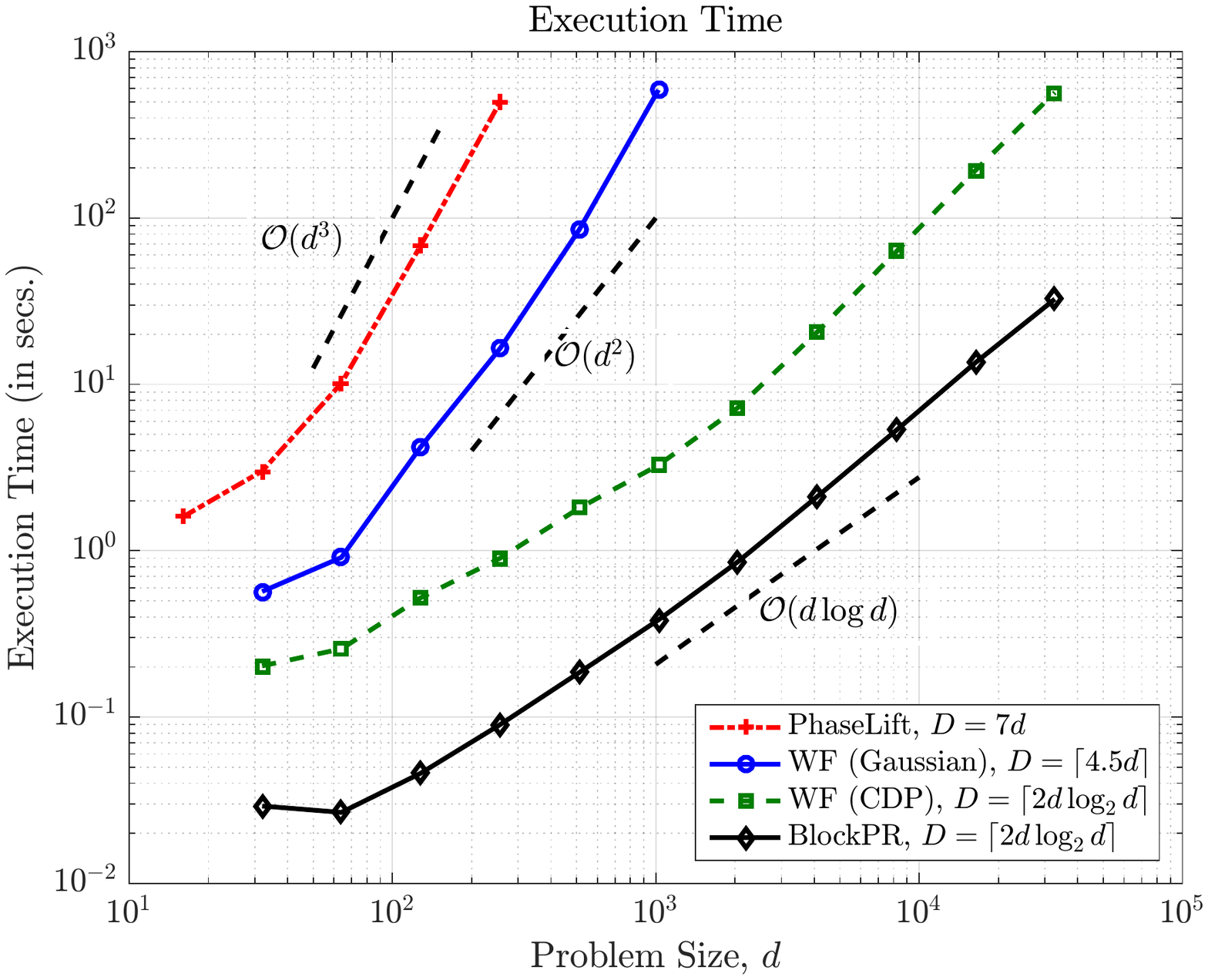}
\caption{Execution time -- Phase Retrieval from $D=\lceil 2d
  \log_2 d\rceil$ measurements.}
\label{fig:exectime}
\end{subfigure}
\caption{Computational Efficiency of the {\em BlockPR} Phase 
Retrieval Algorithm (WF denotes the {\em Wirtinger Flow}
algorithm; parentheses are used to denote the type of
measurements used)}
\label{fig:efficiency}
\end{figure}

We start by presenting numerical simulations demonstrating the
efficiency of the {\em BlockPR} algorithm. In particular, we plot the
execution time for solving the phase retrieval problem (averaged over
100 trials) from perfect (noiseless) measurements in Figure
\ref{fig:efficiency}.  Simulations were performed on a laptop computer
with an Intel$^\circledR$ Core\texttrademark i3-2350M processor, 6GB RAM
and Matlab R2015a. For comparison, we also plot execution times for the
{\em Gerchberg--Saxton} alternating projection algorithm, semidefinite
programming-based {\em PhaseLift} algorithm\footnote{The PhaseLift
    algorithm was implemented as a trace-regularized least-squares
    problem using CVX \cite{cvx,gb08}, a package for specifying and
solving convex programs in Matlab.} and the recently introduced {\em
Wirtinger Flow} algorithm. Simulation results with {\em PhaseLift} and
the {\em Gerchberg--Saxton} alternating projection algorithm use random
complex Gaussian measurements. We note that even though alternating
projection algorithms are known to empirically work with certain
FFT-time measurement constructions, recovery guarantees are only
available for random Gaussian measurements (see for example
\cite{AltProjGaussianMeas}).  Simulation results for the {\em Wirtinger
Flow} algorithm are generated using either random complex Gaussian
measurements or coded diffraction patterns(CDPs). We remark here that
the motivating application typically determines the type of measurements
used.  For example, random Gaussian measurements or coded diffraction
patterns -- such as those employed by {\em PhaseLift} and {\em Wirtinger
Flow} -- are not directly applicable in local correlation-based
applications, which are of primary interest in this discussion.
Nevertheless, Figure \ref{fig:efficiency} provides a useful comparison
for evaluating the efficiency of different phase retrieval algorithms. 
Figure \ref{fig:highsnr_runtime} plots the execution time for
solving the phase retrieval problem using $7d$ measurements. We
observe that {\em PhaseLift} and alternating projections have
computational complexities which scale cubically with the problem
dimension, thereby limiting their applicability to small-scale
problems. On the other hand, both {\em BlockPR} and {\em
Wirtinger Flow} (with coded diffraction pattern measurements)
have essentially FFT-time computational complexities, with 
{\em BlockPR} observed to have a smaller constant than {\em
Wirtinger Flow}. Similar speedup is observed in Figure
\ref{fig:exectime}, where the number of measurements used by {\em
BlockPR} and {\em Wirtinger Flow} grows log-linearly in the
problem size $d$. For reference, execution times for phase
retrieval from $\mathcal O(d)$ random complex Gaussian magnitude
measurements using {\em PhaseLift} and {\em Wirtinger Flow} are
also provided. These plots confirm the significant speedup
offered by the proposed {\em BlockPR} algorithm over other
comparable methods and corroborate the efficiency of the proposed
computational framework for solving large-scale phase retrieval
problems.

\begin{figure}[hbtp]
\centering
\begin{subfigure}[b]{0.475\textwidth}
\includegraphics[clip=true, trim = 0.7in 2in 0.25in
2.5in,scale=0.45]{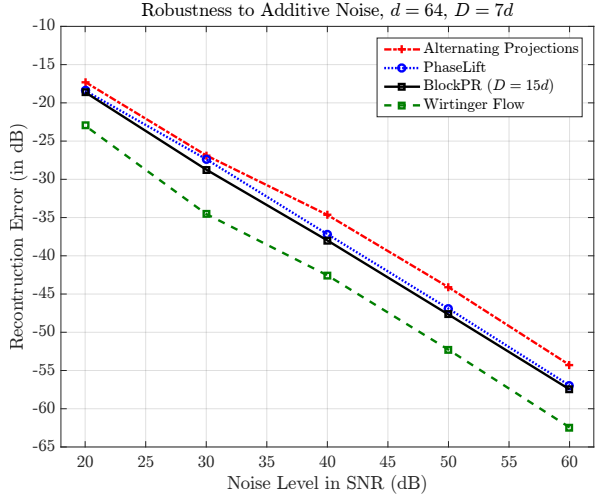}
\caption{Robustness to Additive Noise ($d=64$).}
\label{fig:robustness_d64}
\end{subfigure}
\hfill
\begin{subfigure}[b]{0.475\textwidth}
\includegraphics[clip=true, trim = 0.7in 2in 0.25in
2.5in,scale=.45]{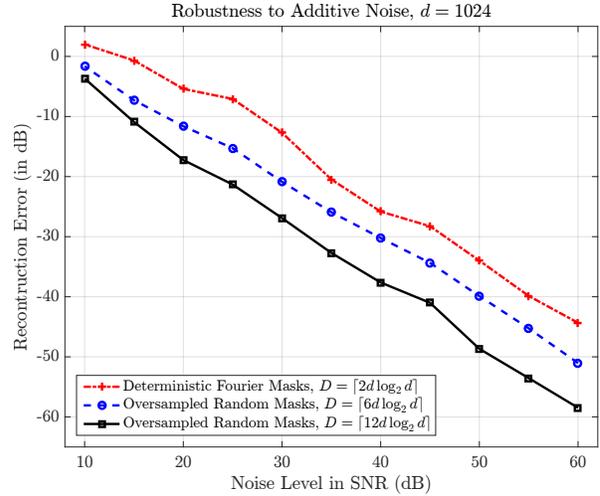}
\caption{Robustness to Additive Noise ($d=1024$).}
\label{fig:robustness_d1024}
\end{subfigure}
\caption{Robustness of the Proposed {\em BlockPR} Phase 
Retrieval Algorithm to Additive Noise}
\label{fig:robustness}
\end{figure}
We next demonstrate robustness of the proposed framework to
additive noise. Figure \ref{fig:robustness_d64} plots the
reconstruction error in recovering a $d=64$ complex random
Gaussian signal at different SNRs, with each data point computed
as the average of 100 trials.\footnote{A few iterations of the
  {\em Gerchberg--Saxton} alternating projection algorithm were
used to post-process the reconstructions.} We include
reconstruction results using {\em Gerchberg--Saxton}, {\em
PhaseLift} and {\em Wirtinger Flow} algorithms for comparison. In
all cases, the algorithms recover the unknown signal $\x$ upto
(or slightly better than) the added noise level. As opposed to
the other algorithms in Figure \ref{fig:robustness_d64}, {\em
BlockPR} uses {\em local} measurements as prescribed in
(\ref{equ:MeasDef}), with each measurement only providing
information about a corresponding {\em local} region of the
underlying signal. Hence, we expect {\em BlockPR} to require more
measurements than the other algorithms in Figure
\ref{fig:robustness_d64} to achieve the same noise robustness.
Indeed, this is confirmed in the figure with {\em BlockPR}
requiring roughly twice the number of measurements to achieve the
same noise robustness. We note, however, that the number of
additional measurements is typically a small constant
irrespective of the problem size. Similar results are provided in
Figure \ref{fig:robustness_d1024} for a larger problem size
($d=1024$). Until now, the {\em deterministic} Fourier-based
measurement masks of \eqref{equ:MeasDef} have been utilized to
generate simulation results. In Figure
\ref{fig:robustness_d1024}, we additionally provide robustness
results using oversampled random local masks. Specifically, we
choose the non-zero entries of each of the $\ell$ measurement
masks $\m_\ell \in \mathbbm C^d$ in (\ref{equ:SubMatDef}) to 
be i.i.d. complex random Gaussian entries for $\ell = 1,2,
\cdots, \gamma\cdot(2\delta-1)$, where $\gamma\geq 1$ is an 
oversampling factor. The dashed line in Figure
\ref{fig:robustness_d1024} was generated using $\gamma=1.5$,
while the solid line was generated using $\gamma=2$. In either 
case, the block length $\delta$ was chosen so as to achieve a
total of $D=\lceil 6d \log_2d\rceil$ and $D=\lceil 12d \log_2
d\rceil$ measurements respectively. This plot confirms the
graceful degradation of reconstruction error with added noise.
Moreover, the block length $\delta$ and oversampling factor
$\gamma$ may be suitably chosen to achieve desired noise 
robustness.
\begin{figure}[htbp]
\centering
\includegraphics[clip=true, trim = 0.75in 2in 0.25in
2.5in,scale=.45]{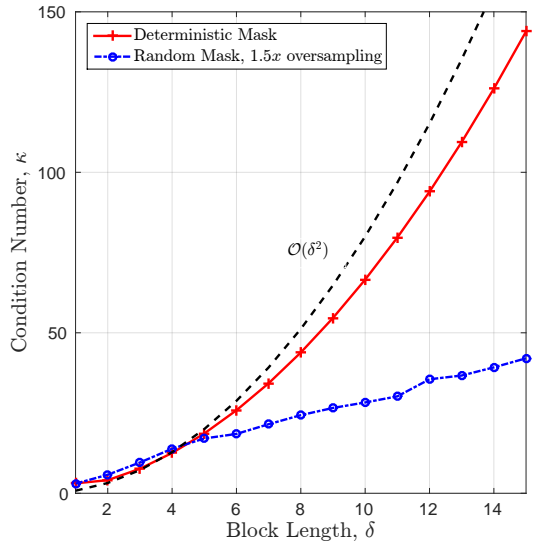}
\caption{Well Conditioned Measurements -- Condition Number as a Function of the Block Length $\delta$}
\label{fig:condno}
\end{figure}

Finally, Figure \ref{fig:condno} plots the condition number of
the system matrix used to solve for the phase differences (matrix
$M'$ in \S \ref{sec:prelims}). The figure plots the condition
number as a function of the block length $\delta$ for
$d=64$.\footnote{ The condition number is independent of the
problem dimension $d$ and depends only on the block length
$\delta$.} It confirms that the condition number scales as a
small multiple of $\delta^2$. The figure also includes a plot of
the condition number when using random masks at an oversampling
factor of $1.5$. Empirical simulations suggest that the use of
oversampling can lead to essentially-linear growth in the
condition number $\kappa(M^\prime)$ with block length $\delta$.
Analyzing the performance of such an oversampled measurement
construction and explicitly writing down associated condition
number bounds would be an interesting avenue for future research.

%
\begin{figure}[hbtp]
\centering
\begin{subfigure}[b]{0.475\textwidth}
\includegraphics[clip=true, trim = 0.65in 2in 0.25in
2.45in,scale=.45]{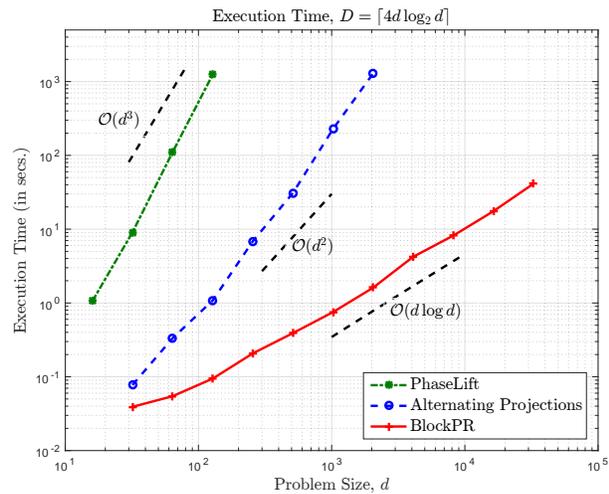}
\caption{Execution Time vs. Problem Size (Local Correlation Measurements).}
\label{fig:exec-local}
\end{subfigure} 
\hfill
\begin{subfigure}[b]{0.475\textwidth}
\includegraphics[clip=true, trim = 0.65in 2in 0.25in
2.5in,scale=.45]{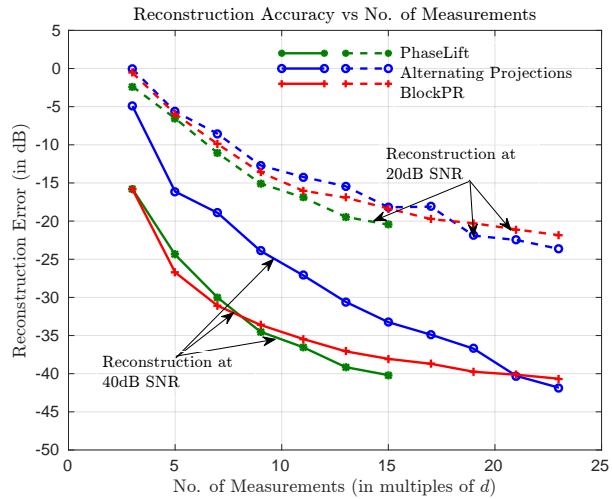}
\caption{Reconstruction Error vs. Number of Measurements ($d=64$,
i.i.d. Gaussian noise).}
\label{fig:transition-noise}
\end{subfigure} \vspace{0.1in} \\
\begin{subfigure}[b]{0.475\textwidth}
\includegraphics[clip=true, trim = 0.75in 2in 0.25in
2.5in,scale=0.45]{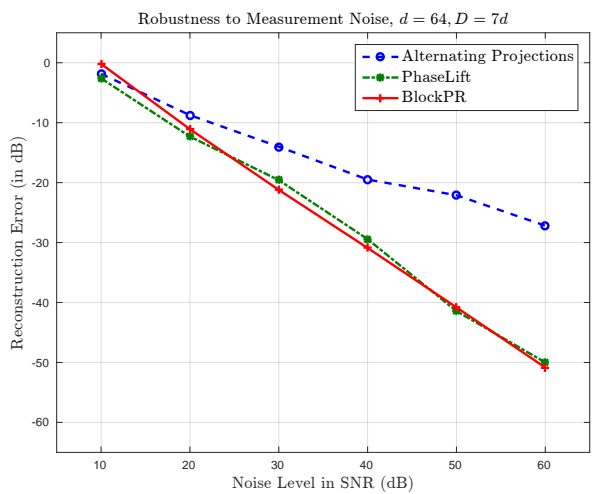}
\caption{Robustness to Additive Noise ($d=64, D=7d$, (Fourier-like) 
    Local Correlation Measurements)}
\label{fig:robustness-local-7d}
\end{subfigure}
\hfill
\begin{subfigure}[b]{0.475\textwidth}
\includegraphics[clip=true, trim = 0.75in 2in 0.25in
2.5in,scale=0.45]{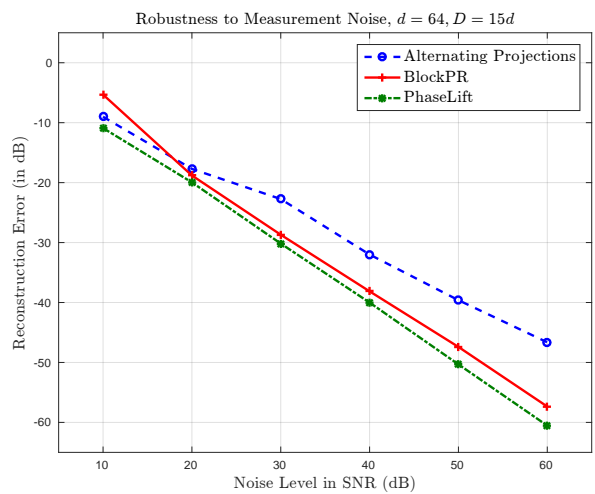}
\caption{Robustness to Additive Noise ($d=64, D=15d$, (Fourier-like) 
    Local Correlation Measurements)}
\label{fig:robustness-local-15d}
\end{subfigure}
\caption{Evaluation of Phase Retrieval Algorithms for {\em Local}
Correlation Measurements}
\label{fig:local-performance}
\end{figure}
We conclude this section by presenting representative numerical results
highlighting the superior performance of the proposed {\em BlockPR}
algorithm for {\em local} measurements. Specifically, we use the
windowed Fourier-like local correlation measurements described in
Section \ref{sec:CondNumBounds} with the {\em BlockPR}, {\em PhaseLift} and {\em
Gerchberg-Saxton} alternating projection algorithms. We note that there
are no theoretical recovery guarantees with {\em PhaseLift} and the {\em
Gerchberg-Saxton} algorithm for such measurements; however, empirical
simulations suggest that for sufficiently large numbers of measurements,
both algorithms recover the unknown signal. Indeed, this is illustrated
in Figure \ref{fig:transition-noise}, where the reconstruction error in
recovering a $d=64$ length complex vector is plotted as a function of
the number of measurements required. For global (random Gaussian)
measurements, it is known that {\em PhaseLift} recovers signals upto the
noise level with about $6d$ measurements (see, for example
\cite{candes2013phaselift}). Yet, for the local measurements studied
here, significantly larger numbers of measurements are necessary for
successful phase retrieval by {\em PhaseLift} (and indeed the {\em
Gerchberg-Saxton} alternating projection algorithm). Figure
\ref{fig:transition-noise} shows that the proposed {\em BlockPR}
algorithm compares well with {\em PhaseLift} and almost always
outperforms alternating projections. We note that the {\em PhaseLift}
error plots are only shown for upto $D=15d$ measurements due to memory
issues in implementing the algorithm for larger $D$ with CVX on a laptop
computer. 

The performance of the proposed {\em BlockPR} algorithm is particularly
noteworthy since it has an essentially FFT-time computational cost as
illustrated in Figure \ref{fig:exec-local}. Unsurprisingly, the
computational cost of {\em PhaseLift} scales as $\mathcal O(d^3)$. Note
that we had to use the Matlab software package TFOCS
\cite{TFOCS-paper,TFOCS-web} to implement {\em PhaseLift} for large
problem sizes. TFOCS implements fast first-order conic solvers which
trade off some accuracy for efficiency gains, thereby enabling the
solution of large problems.  We also note that the computational cost of
the {\em Gerchberg-Saxton} alternating projections algorithm scales
quadratically despite the use of FFT-based linear system solvers in each
iteration. Numerical experiments suggest that the number of alternating
projection iterations grows with the problem size, thereby resulting in
the quadratic scaling.

For completeness, Figures \ref{fig:robustness-local-7d} and
\ref{fig:robustness-local-15d} plot the reconstruction error versus
added noise level using $D=7d$ and $D=15d$ measurements respectively for
all three algorithms using local correlation measurements. In each plot,
we recover a $d=64$ length complex vector from i.i.d Gaussian noise
corrupted phaseless measurements. With $D=15d$ measurements, both {\em
BlockPR} and {\em PhaseLift} recover the unknown signal to the level of
noise, while for $D=7d$ measurements, there is a roughly 10dB error in
the reconstructed signal. Note the increased number of measurements
necessary for {\em PhaseLift} to achieve the same error performance as
in Figure \ref{fig:robustness_d64} -- this highlights the challenging
nature of reconstructing signals from {\em local} phaseless
measurements.

\section{Extension: Sublinear-Time Phase Retrieval for Compressible Signals}
\label{sec:PRsparse}

In this section we briefly focus on the compressive phase retrieval
setting, (see, e.g.,
\cite{ohlsson2012cprl,SCHNITNER!,li2013sparse,wang2014phase,eldar2014phase,GESPAR2014Eldar}),
where one aims to approximate a sparse or compressible ${\bf x} \in
\mathbbm{C}^d$ using fewer magnitude measurements than required for the
recovery of general ${\bf x}$.  It is known that robust compressive
phase retrieval for $s$-sparse vectors is possible using only
$\mathcal{O}(s\log(d/s))$ magnitude measurements
\cite{eldar2014phase,Iwen2014robust}.  In this section we prove that it
is also possible to recover $s$-sparse vectors ${\bf x} \in
\mathbbm{C}^d$ up to an unknown phase factor in only $\mathcal{O}(s
\log^6 d)$-time using $\mathcal{O}(s \log^5 d)$ magnitude measurements.
Thus, we establish the first known nearly runtime-optimal (i.e.,
essentially linear-time in $s$) compressive phase retrieval recovery
result.  In particular, we prove the following theorem.

\begin{thm}
There exists a deterministic algorithm $\mathcal{A}: \mathbbm{R}^D \rightarrow \mathbbm{C}^d$ for which the following holds:  Let $\epsilon \in (0,1]$, ${\bf x} \in \mathbbm{C}^{d}$ with $d$ sufficiently large, and $s \in [d]$.  Then, one can select a random measurement matrix $\tilde{M} \in \mathbbm{C}^{D \times d}$ such that 
\begin{equation}
\min_{\theta \in [0, 2 \pi]} \left\| \mathbbm{e}^{\mathbbm{i} \theta} {\bf x}  - \mathcal{A} \left( |  \tilde{M} \x |^2 \right) \right\|_2~\leq~ \left\| {\bf x} - {\bf x}^{\rm~ opt}_s \right\|_2 + \frac{22 \epsilon \left\| {\bf x} - {\bf x}^{\rm opt}_{(s/\epsilon)} \right\|_1}{\sqrt{s}}
\label{thm:ArbSparseRecovery}
\end{equation}
is true with probability at least $1-\frac{1}{C \cdot \ln^2(d) \cdot \ln^3 \left( \ln d \right)}$.\footnote{Here $C \in \mathbbm{R}^+$  is a fixed absolute constant.}
Here $D$ can be chosen to be $\mathcal{O} \big( \frac{s}{\epsilon} \cdot \ln^3( \frac{s}{\epsilon} ) \cdot \ln^3 \left( \ln \frac{s}{\epsilon} \right)$ $\cdot \ln d \big)$.  Furthermore, the algorithm will run in $\mathcal{O} \left( \frac{s}{\epsilon} \cdot \ln^4( \frac{s}{\epsilon} ) \cdot \ln^3 \left( \ln \frac{s}{\epsilon} \right) \cdot \ln d \right)$-time in that case.\footnote{For the sake of simplicity, we assume $s = \Omega(\log d)$ when stating the measurement and runtime bounds above.}
\label{thm:RecoverArbSparseVecs}
\end{thm}

We prove Theorem~\ref{thm:RecoverArbSparseVecs} by following the generic compressive phase retrieval recipe presented in \cite{Iwen2014robust}.  Let $C \in \mathbbm{C}^{m \times d}$ be any compressive sensing matrix with an associated sparse approximation algorithm $\Delta: \mathbbm{C}^m \rightarrow \mathbbm{C}^d$ (see, e.g., \cite{CS1,NearOpt,CS2,COSAMP,HardThreshforCS,ROMP,ROMPstable}), and let $P \in \mathbbm{C}^{D \times m}$ be any phase retrieval matrix with an associated recovery algorithm $\Phi:  \mathbbm{R}^D \rightarrow \mathbbm{C}^m$.  Then, $\Delta \circ \Phi:  \mathbbm{R}^{D} \rightarrow \mathbbm{C}^d$ will approximately recover compressible vectors ${\bf x} \in \mathbbm{C}^d$ up to an unknown phase factor when provided with the magnitude measurements $|PC\x|$.  That is, one may first use $\Phi$ to recover $\mathbbm{e}^{\mathbbm{i} \phi} (C{\bf x}) ~=~C ( \mathbbm{e}^{\mathbbm{i} \phi}  {\bf x})$ for some unknown $\phi \in [0, 2\pi]$ from $|PC\x|$, and then use $\Delta$ to recover $\mathbbm{e}^{\mathbbm{i} \phi}  {\bf x}$ from $C ( \mathbbm{e}^{\mathbbm{i} \phi}  {\bf x} )$.  If both $\Phi$ and $\Delta$ are efficient, the result will be an efficient sparse phase retrieval method.

Herein we will utilize Algorithm~\ref{alg:phaseRetrieval} as our phase retrieval method.  Note that it's runtime is only $\mathcal{O}(m \log^4 m)$, making it optimal up to log factors (recall Theorem~\ref{thm:RecoverArbVecs}).  For the compressive sensing method we will utilize the following algorithmic result from \cite{Iwen2014compressed}.
\begin{thm}
Let $\epsilon \in (0,1]$, $\sigma \in [2/3,1)$, ${\bf x} \in \mathbbm{C}^{d}$, and $s \in [d]$.  With probability at least $\sigma$ the deterministic compressive sensing algorithm from \cite{Iwen2014compressed} will output a vector ${\bf z} \in \mathbbm{C}^d$ satisfying
\begin{equation}
\left\| {\bf x} - {\bf z} \right\|_2 ~\leq~ \left\| {\bf x} - {\bf x}^{\rm~ opt}_s \right\|_2 + \frac{22 \epsilon \left\| {\bf x} - {\bf x}^{\rm opt}_{(s/\epsilon)} \right\|_1}{\sqrt{s}}
\label{equ:ApproxError}
\end{equation}
when executed with random linear input measurements $\mathcal{M}{\bf x} \in \mathbbm{C}^m$.  Here $m = \mathcal{O} \left(\frac{s}{\epsilon} \cdot \ln \left( \frac{s / \epsilon }{1-\sigma} \right) \ln d \right)$ suffices.  The required runtime of the algorithm is $\mathcal{O} \left(\frac{s}{\epsilon} \cdot \ln \left( \frac{s / \epsilon }{1-\sigma} \right) \ln \left( \frac{d}{1-\sigma} \right) \right)$ in this case.\footnote{For the sake of simplicity, we assume $s = \Omega(\log d)$ when stating the measurement and runtime bounds above.}  
\label{appthm:NewRecoverAlg}
\end{thm}
Theorem~\ref{thm:RecoverArbSparseVecs} now follows easily from Theorem~\ref{thm:RecoverArbVecs} with $\n = {\bf 0}$, and Theorem~\ref{appthm:NewRecoverAlg}.

\section{Future Work}
\label{sec:FutureWork}

This paper provides the first known global robust recovery guarantees for ptychographic phase retrieval problems, as well a new numerical lifting+angular synchronization solution approach which is applicable with correlation-type measurements resulting from a wide class of locally supported masks $\m_\omega$.  In doing so, it opens up many additional avenues of research.  Examples include, e.g., questions associated with the use of different boundary conditions for $\x$.  Herein, $\x$ is assumed to be periodic, which leads to the lifted linear system matrix $M'$ in \eqref{equ:LinProb} being block circulant.  It would be interesting to derive condition number bounds, as well as specialized numerical solution techniques, for the lifted linear system matrices that result from alternate boundary condition assumptions.  More generally, it would be interesting to derive theoretical condition number bounds for the types of lifted linear systems that arise from more general classes of masks $\m_\omega$ under various sets of boundary condition assumptions for $\x$. 

It would also be interesting to explore alternate approaches, besides Algorithm~\ref{alg:angSync}, for solving the angular synchronization problem via local phase differences that appears in the second stage of our proposed phase retrieval approach.  In particular, it would be interesting to develop theoretical error bounds for other solution methods (e.g., \cite{ang_sync_singer}) in the presence of noise for these types of local angular synchronization problems.  Finally, it would also be instructive to carefully consider the optimality (or lack thereof) of the robustness guarantees derived herein for $m$-flat signals (i.e., Theorem~\ref{thm:RecoverFlatVecs}) under various noise models.

\section*{Acknowledgements}

The authors would like to thank the anonymous reviewers for many helpful suggestions and interesting ideas regarding possible directions for future work.

\bibliographystyle{siam}
\bibliography{FastPR}

\end{document}